\documentclass[11pt]{article}
\usepackage[T1]{fontenc}
\usepackage[utf8]{inputenc}
\usepackage{lmodern}
\pdfoutput=1
\usepackage{amsmath,amssymb,amsthm,mathtools}
\usepackage{microtype}
\usepackage{enumitem}
\usepackage[a4paper,margin=1in]{geometry}
\usepackage{cite}
\usepackage[hidelinks]{hyperref}
\usepackage{mathrsfs}
\usepackage{bm}
\usepackage{graphicx}
\usepackage{placeins}
\usepackage{float}
\graphicspath{{./}}

\newcommand{\paperfigure}[2][]{%
  \IfFileExists{#2}{\includegraphics[#1]{#2}}{%
    \fbox{\parbox{0.90\textwidth}{\centering\small
      Figure file \texttt{\detokenize{#2}} not found.\\
      Use the complete source package to reproduce the figure.}}%
  }%
}

\newcommand{\C}{\mathbb C}
\newcommand{\R}{\mathbb R}
\newcommand{\Z}{\mathbb Z}
\newcommand{\N}{\mathbb N}
\newcommand{\Q}{\mathbb Q}
\newcommand{\Pj}{\mathbb P}
\newcommand{\Nalg}{\mathcal N}

\newcommand{\Id}{I}
\newcommand{\one}{\mathbf 1}
\newcommand{\wrp}{\wr}

\newcommand{\norm}[1]{\lVert#1\rVert}

\newcommand{\smin}{s_{\min}}

\newcommand{\Tr}{\operatorname{Tr}}

\newcommand{\Stab}{\operatorname{Stab}}
\newcommand{\diag}{\operatorname{diag}}

\theoremstyle{plain}
\newtheorem{theorem}{Theorem}[section]
\newtheorem{proposition}[theorem]{Proposition}
\newtheorem{lemma}[theorem]{Lemma}
\newtheorem{corollary}[theorem]{Corollary}

\newtheorem{definition}[theorem]{Definition}
\theoremstyle{definition}
\newtheorem{conjecture}[theorem]{Conjecture}
\newtheorem{question}[theorem]{Question}

\newtheorem{introdefinition}{Definition}
\newtheorem{example}[theorem]{Example}

\theoremstyle{remark}
\newtheorem{remark}[theorem]{Remark}

\newenvironment{keyword}{\par\vspace{0.5em}\noindent\textit{Keywords:}\ }{\par\vspace{0.5em}}

\title{Spectral Dynamics and Fuglede--Kadison Determinants for Wreath Products}
\author{Penghui Wang \ \ Xumin Wang\footnote{Corresponding author.} \ \ Rongwei Yang}
\date{}

\begin{document}
\maketitle

\begin{abstract}
We investigate connections between operator spectra, complex dynamics, and determinants for wreath product groups. For a nontrivial finite abelian group $G$, the natural action of $G\wr\mathbb Z$ on a rooted tree produces a family of operators with a recursive structure. Applying Schur complements to this family yields a rational map on complex projective space. We determine the projective spectrum—the parameter values for which the operator pencil fails to be invertible—as the union of  the Julia set of this map and an explicitly described subset of its Fatou set. We also prove that the associated Koopman representation and the left regular representation are weakly equivalent and therefore give the same projective spectrum.

We then compute Fuglede–Kadison determinants for a family of operators on more general wreath products. These formulas yield injective operators with zero determinant for wreath products $K\wr\mathbb Z$, where $K$ is a discrete group admitting a nontrivial projection $\mathfrak{e}\in \C[K]$. This answers negatively a question of Li and Thom about whether injectivity implies positive determinant for operators defined by matrices over complex group rings of amenable groups. For $G_{tf}\rtimes_\sigma\mathbb Z$, where $G_{tf}$ is any countable torsion-free abelian group and $\sigma\in\operatorname{Aut}(G_{tf})$, we prove that injectivity does imply positive determinant in every matrix size.
\end{abstract}

\begin{keyword}
lamplighter groups, Wreath Products, Koopman representation, left regular representation, projective spectrum, spectral dynamics, Fuglede-Kadison determinant
\end{keyword}

\section{Introduction}\label{s1}

Let $\mathscr{G}$ be a countable discrete group acting on a standard probability space $(X,\mu)$ by measure-preserving transformations. The associated \emph{Koopman representation} is the unitary representation $\pi$ on $L^2(X,\mu)$ defined by
\begin{equation}\label{e1}
	(\pi(g)f)(x)=f(g^{-1}x).
\end{equation}
Koopman's formulation of Hamiltonian dynamics \cite{Koopman} and von Neumann's mean ergodic theorem \cite{vonNeumann} are early examples of the operator-theoretic approach to ergodic theory. In this setting, ergodicity means that the only invariant functions are the constants, while weak mixing is equivalent to the absence of nonzero finite-dimensional invariant subspaces in
\[
L_0^2(X,\mu)=\left\{f\in L^2(X,\mu)\mid\int_X f\,d\mu=0\right\}.
\]
The Halmos--von Neumann theorem states that ergodic transformations with discrete spectrum are rotations on compact abelian groups \cite{HalmosvonNeumann}. Further relations with representation theory, measurable cocycles, and rigidity were developed by Mackey and Zimmer \cite{Mackey,Zimmer}. Koopman operators also occur in nonlinear model reduction and the study of fluid flows \cite{Mezic,Rowley}.

The left regular representation of $\mathscr{G}$ is
\begin{equation}\label{e2}
	(\lambda_\mathscr{G}(g)\xi)(h)=\xi(g^{-1}h),\ \ \xi\in\ell^2(\mathscr{G}).
\end{equation}
It is the Koopman representation associated with the left-translation action of $\mathscr{G}$ on itself equipped with counting measure. Thus the regular representation is a special case of the Koopman representation, whereas a general Koopman representation depends on the chosen action. The constant function $1$ is invariant under $\pi_\mu(\mathscr{G})$, whereas $\lambda_\mathscr{G}$ has no nonzero invariant vector when $\mathscr{G}$ is infinite. Weak containment provides a useful relation between representations of infinite groups. We write $\rho_1\prec\rho_2$ if
\[
\|\rho_1(a)\|\le \|\rho_2(a)\|,\ \ a\in\C[\mathscr{G}].
\]
If both $\rho_1\prec\rho_2$ and $\rho_2\prec\rho_1$, the representations are \emph{weakly equivalent}, written as $\rho_1\sim_w\rho_2$. Equivalently, the corresponding integrated representations of $C^*(\mathscr{G})$ have the same kernel \cite{Fell}. Hulanicki characterized amenability in terms of weak containment in the regular representation \cite{Hulanicki}, and Kesten obtained the corresponding spectral-radius criterion for symmetric random walks \cite{Kesten}. Dudko--Grigorchuk studied weak containment for Koopman, groupoid, and orbit representations. In particular, the Koopman representation of an essentially free measure-preserving action of an amenable group is weakly equivalent to the regular representation \cite{DudkoGrig}.

For a tuple $A=(A_0,\ldots,A_m)$ of elements in a Banach algebra, define the multiparameter pencil
\[
A(z)=z_0A_0+\cdots+z_mA_m.
\]
Its \emph{projective spectrum} is defined by
\[
p(A)=\left\{[z_0:\cdots:z_m]\in\Pj^m\mid A(z)\text{ is not invertible}\right\},
\]
and it connects several complex variables, cohomology, Hermitian geometry, operator theory, operator algebras, representation theory, and complex dynamics \cite{CadeYang,Chagouel2015,DouglasYang,Feng2026,GengLiuWang2024,GY,GoldbergYang, HeWangYang, HeYang, MaoQiaoWang2017,YangProj}. Monograph \cite{YangBook} provides a good account of these developments. Given a finitely generated group $H=\langle g_1, ..., g_m\rangle$ and a unitary representation $\rho$ of $H$, the projective spectrum of the pencil 
\[A_\rho(z)=z_0I+z_1\rho(g_1)+\cdots +z_m\rho(g_m)\] is known to encode structural information about $H$ and $\rho$. In particular, we have
\[
\rho_1\prec\rho_2
\Longrightarrow
p(A_{\rho_1})\subseteq p(A_{\rho_2}),
\]
and consequently,
\[
\rho_1\sim_w\rho_2
\Longrightarrow
p(A_{\rho_1})=p(A_{\rho_2}).
\]
The same assertion holds for matrix-valued pencils \cite[Sec.~2]{GY}.

We use the interval model in \cite[Chap.~10]{YangBook} to describe group actions on a $q$-nary tree. Let $X$ be a finite alphabet with $|X|=q\ge2$, equipped with the usual order $X=\{0,1,\ldots,q-1\}$. Apart from the usual $q$-adic endpoints, an infinite word $x_1x_2\cdots\in X^{\mathbb N}$ corresponds to
\[
\sum_{k\ge1}\frac{x_k}{q^k}\in[0,1).
\]
A finite word $w=x_1\cdots x_n,\ n\geq 1,$ determines a $q$-adic subinterval of length $q^{-n}$. If a group $\mathscr{G}$ acts on $[0, 1)$ by permuting these $q$-adic subintervals, then the associated Koopman representation of $\mathscr{G}$ on $\mathscr H:=L^2[0,1)$ is defined by
\[
(\pi(g)f)(x)=f(g^{-1}x),\ \ \ g\in \mathscr{G}.
\]

The partition of $[0,1)$ into $q$ subintervals of length $q^{-1}$ yields a unitary map $W: \mathscr H\to \mathscr H^{\oplus q}$ defined by
\[
(Wf)_j(x)=q^{-1/2}f\!\left(\frac{j+x}{q}\right),\ \
0\le j\le q-1.
\]
We say that the representation $\pi$ is self-similar if the operators $W\pi(g)W^*,\ g\in \mathscr{G}$, are given by  $q\times q$ block matrices whose nonzero blocks are all of the form $\pi(h)$ for some $h\in \mathscr{G}$. 
For $n\ge0$, denote by $X^n$ the set of $q$-adic subintervals of length $q^{-n}$ and
let $\mathscr H_n\subset \mathscr H$ be the subspace of functions constant on each subinterval in $X^n$. Then $\mathscr H_n\cong\ell^2(X^n)$, the sequence $\mathscr H_0\subset \mathscr H_1\subset\cdots$ has dense union in $\mathscr H$, and the restriction of the Koopman representation $\pi$ to $\mathscr H_n$, denoted by $\pi_n$, is a permutation representation of $\mathscr{G}$ on $X^n$.

Projective spectra connect self-similar group representations with complex dynamics. Applying Schur complement argument to the pencil $A_{\pi_n}(z)$ often yields a rational map
\[
\widehat F:\Pj^m\rightarrow\Pj^m.
\]
Its indeterminacy set consists of the points where all homogeneous coordinate functions vanish. Away from the indeterminacy sets of the iterates, a point belongs to the Fatou set ${\mathcal F}(\widehat F)$ if it has a neighborhood on which the iterates form a normal family, the Julia set ${\mathcal J}(\widehat F)$ is the complement of the Fatou set. 
In Subsection~\ref{ss1}, we define our convention following the spectral-dynamics framework of \cite{GrigNek,GoldbergYang,YangBook}.

The finite-dimensional representations $\pi_n$ provide a spectral approximation of $\pi$ in the sense that, if the pencil $\pi(h)$ is selfadjoint, where $h\in \C[\mathscr{G}]$, then its classical spectrum \[\sigma(\pi(h))=\overline{\cup_{n=1}^\infty \sigma (\pi_n(h))},\] see \cite{GNS}. Quasi-regular representations give another description through random walks on orbit Schreier graphs. For a finite symmetric generating set $\mathcal S$ of group $\mathscr{G}$, define the Hecke-type operator (also called the Markov operator)
\[
M_\pi=|\mathcal S|^{-1}\sum_{s\in\mathcal S}\pi(s).
\]
If $\pi_x$ is the orbit representation at $x\in X^{\mathbb N}$, then the inner product $\langle M_{\pi_x}^n\delta_x,\delta_x\rangle$ gives the $n$-step return probability at $x$ \cite{GNS}. In summary, the restrictions $\pi_n$, the orbit representations, and the rational map $F$ reveal three different aspects of the spectral property of $\mathscr{G}$ and $\pi$. Describing their connection is an appealing undertaking.

Bartholdi--Grigorchuk computed the spectra of Hecke-type operators for fractal groups via finite-dimensional approximation and obtained examples related to quadratic Julia sets and Cantor spectra \cite{BartholdiGrig}. Grigorchuk--Nekrashevych showed that Schur reduction translates spectral questions for Markov-type elements in self-similar $C^*$-algebras into questions about invariant sets of multivariable rational maps \cite{GrigNek}. A converse relation also holds: Nekrashevych proved that for a postcritically finite rational map of degree at least two, the limit dynamical system of its iterated monodromy group is topologically conjugate to the map on its Julia set \cite[Theorem~6.4.4]{NekBook}. Thus self-similarity encapsulates both Julia set and its dynamics.

Schreier graphs provide algebraic models for the Sierpi\'nski gasket and the Basilica Julia set \cite{GNS}. For the Sierpi\'nski gasket, Fukushima--Shima used spectral decimation to compute the eigenvalues of the Laplacian and their asymptotic distribution \cite{FukushimaShima}. For the Basilica group, Bartholdi--Vir\'ag used self-similarity and random walk estimates to prove amenability \cite{BartholdiVirag}. These examples illustrate the close connections among recursive group actions, spectral theory, geometry, and probability.

Wreath products form a particularly explicit class of self-similar groups. For a finite abelian group $G$, Silva--Steinberg and Kambites--Silva--Steinberg realize the wreath product $G\wr\mathbb{Z}$ as a self-similar group acting on the rooted tree $G^*$, or equivalently on its boundary $G^{\mathbb{N}}$, via the Cayley automaton \cite{SilvaSteinberg,KSS}. We use this realization to obtain recursive matrix representations on $\ell^2(G^n)$. Grigorchuk--\.{Z}uk considered the lamplighter group, i.e., the case $G=\mathbb{Z}_2$, and obtained a purely atomic spectral measure for the natural Markov operator \cite{GZ}. Kambites--Silva--Steinberg extended the corresponding calculations to general finite groups $\mathscr{G}$ \cite{KSS}. For $\mathscr{G}\wr H$ with $H$ finitely generated, Lehner--Neuhauser--Woess related the spectral measure of the symmetric switch--walk--switch operator to an absorbing random walk on a Bernoulli percolation cluster with parameter $|\mathscr{G}|^{-1}$, almost sure finiteness of the clusters yields pure point spectrum \cite{LNW}. Spectral computations of the lamplighter group also appear in the construction by Grigorchuk--Linnell--Schick--\.{Z}uk of a closed manifold with an $L^2$-Betti number equal to $1/3$ \cite{GLSZ}. These results provide explicit operator computations for wreath products in recursive representations, random walks, and $L^2$-invariants.

For self-similar pencils, Schur complement yields a rational map whose Julia set can be compared with the projective spectrum. This relation is referred to as \emph{spectral dynamics} in \cite{GoldbergYang}. For the infinite dihedral group
\[
D_\infty=\langle a,t\mid a^2=t^2=1\rangle,
\]
Grigorchuk--Yang computed the projective spectrum of $z_0I+z_1\lambda(a)+z_2\lambda(t)$, proved that the Koopman and regular representations are weakly equivalent, and obtained a rational map that preserves the spectrum \cite[Theorems~1.1,~7.2 and Sec.~8]{GY}. Together with the Fatou--Julia analysis of \cite{ZYL} and \cite[Chap.~10]{YangBook}, these results lead to the following question.

\medskip
\noindent\textbf{Spectral-dynamics question.}
Given a self-similar representation $\pi$ of group $\mathscr{G}$, how much of the projective spectrum $p(A_\pi)$ is described by the Julia set of the rational map $\widehat F$? In particular, can $p(A_\pi)$ contain Fatou points of $\widehat F$, and if so, how can this part of the spectrum be described?

\medskip

Rational maps arising from Schur reduction are often called renormalization maps. We follow the convention for indeterminacy given in Subsection~\ref{ss1}. Goldberg--Yang related the Julia set to the projective spectrum and the extended indeterminacy set in the dihedral case, obtaining partial results and a conjecture for the Grigorchuk group \cite[Theorem~5.11, Corollary~6.5 and Conjecture~6.6]{GoldbergYang}. Yang \cite[Chap. 10]{YangBook} shows that projective spectrum in fact coincides with the Julia set for the $D_\infty$ case. Dang--Grigorchuk--Lyubich interpreted the spectral measures for the Grigorchuk, lamplighter, and Hanoi Schreier operators as slices of spectral currents obtained from algebraic pullbacks. They distinguished an absolutely continuous density of states in the Grigorchuk case from atomic densities in the lamplighter and Hanoi cases, and obtained quantitative approximation results in the atomic cases \cite[Theorem~A]{DangGL}. For the Basilica group, Liu--He recently described subsets common to the Julia set and the projective spectrum and those common to the Fatou set and the resolvent set; the full equality remains conjectural \cite{LiuHe}. Zu--Yang--Lu proved that adjoining the direction $at$ to the dihedral pencil yields a projective spectrum equal to the Julia set of the induced map \cite[Theorem~A]{ZYL}. For a pencil associated with the lamplighter group $\mathbb{Z}_2\wr\mathbb{Z}$, they obtained a Julia-set description and proposed
\[
p(D_\pi)=\phi(L)\cup\mathcal{J}(\widehat F)
\]
in \cite[Theorem~B and Conjecture~4.9]{ZYL}.

For a finite abelian group $G$, the wreath product $G\wr\mathbb{Z}$ admits an explicit self-similar action on $G^{\mathbb{N}}$ and an explicit spectral decomposition in the left regular representation \cite{SilvaSteinberg,KSS,DicksSchick}. A main objective of this paper is to compare the rational dynamics obtained from Schur reduction with the projective spectrum of a non-normal pencil.

Let $G=\{g_1={\bf 1},g_2,\ldots,g_r\}$ be a nontrivial finite abelian group and let $\Gamma_G=G\wr\mathbb{Z}$. We consider the self-similar Koopman representation $\pi$ of $\Gamma_G$ on $\mathscr H=L^2[0,1)$ (constructed originally on the rooted tree in Silva--Steinberg \cite[Sec.~3]{SilvaSteinberg} and Kambites--Silva--Steinberg \cite[Sec.~3, especially (3.6)--(3.9)]{KSS}). Let $t$ be the generator of the $\mathbb{Z}$-factor, acting by shift on the base group $\bigoplus_{\mathbb{Z}}G$. The transformations used here are
\begin{equation}\label{e3}
	s_j=g_j^{-1}t,\ \
	s_1=t.
\end{equation}
Let
\begin{equation}\label{e5}
	e=\frac1r\sum_{g\in G}g,\ \ \ \ 
	C_G=re-\one=\sum_{g\in G\setminus\{1\}}g,
\end{equation}
and define 
\begin{equation}\label{e4}
X=r\pi(t^{-1}e),\ \ S=\pi(C_G),\ \ X_n=r\pi_n(t^{-1}e),\ \ S_n=\pi_n(C_G), \ \ n\ge 1.
\end{equation}
Under the left regular representation of $G$, the element $C_G$ is represented by the $r\times r$ matrix $T_r$ with zero diagonal and all off-diagonal entries equal to $1$. Thus,
\begin{equation}\label{e6}
	S_0=r-1, \ \ \text{and}\ \ \ \ S_n=T_r\otimes I_{r^{n-1}}, \ \ \ \forall n\ge 1.
\end{equation}
These matrices coincide with the operators $S_n$ in \cite[Sec.~5.1]{KSS}. 

For any countable group $\mathscr{G}$, the involution on the group algebra $\mathbb{C}[\mathscr{G}]$ is defined by
\[
\bigg(\sum_{g\in \mathscr{G}}a_ g g\bigg)^*
=\sum_{g\in \mathscr{G}}\overline{a_g}\,g^{-1}.
\]

We consider the pencils
\[
D_{\pi,G}(z)=z_0I+z_1S+z_2X+z_3X^*,
\]
\begin{align}\label{e7}
	D_n(z)&=z_0I_{\mathscr H_n}+z_1S_n+z_2X_n+z_3X_n^*,\ \ z\in\C^4.
\end{align}
Observe that these pencils are in general not selfadjoint or normal. The Schur recursion for $D_n(z)$ yields a rational map $\widehat F$ on $\Pj^3$. 
The set $\mathcal F_\infty\subset\mathcal F(\widehat F)$ is defined in Subsection~\ref{ss1}. The first main result provides an exact description of the projective spectrum.

\begin{theorem}\label{t1}
For every nontrivial finite abelian group $G$, we have
\[
p(D_{\pi,G})
=\mathcal J(\widehat F)\cup \mathcal F_\infty.
\]
\end{theorem}

Theorem~\ref{t1} shows that the Julia set does not capture the entire projective spectrum: some points in the Fatou set also lie in $p(D_{\pi,G})$. The set $\mathcal F_\infty$ is described explicitly in Subsection~\ref{ss1}. For $G=\Z_2$, the expression proposed in \cite[Conjecture~4.9]{ZYL} must therefore include an additional subset of the Fatou set.

We next compare the Koopman representation with the left regular representation. This yields another description of the projective spectrum in terms of an explicit planar region.

\begin{theorem}\label{t2}
For every nontrivial finite abelian group $G$, we have $\pi\sim_w\lambda$ for $\Gamma_G$, and
\[
p(D_{\pi,G})=p(D_{\lambda,G})
=\bigl\{[z]\mid z_0+(r-1)z_1\in E(z_2,z_3)\bigr\}
\cup \mathcal J(\widehat F),
\]
where
\[
E(u,v)=\{-r(uw+v\overline w)\mid w\in\C,\ |w|\le1\}.
\]
\end{theorem}

The set $E(u,v)$ is a closed elliptical region, including the degenerate cases. The weak equivalence in Theorem~\ref{t2} implies equality of the projective spectra, but does not imply unitary equivalence of the representations. Theorem~\ref{t1} describes the parts of the spectrum lying in the Julia and Fatou sets, whereas Theorem~\ref{t2} gives a geometric description of the same spectrum and shows that it remains unchanged when $\pi$ is replaced by $\lambda$.

The Jacobi matrices obtained from the Schur recursion also arise in the left regular representation. This suggests studying Fuglede--Kadison (FK) determinants $\det$ alongside projective spectra. Such determinants already occur in Grigorchuk--Yang's study of the dihedral pencil
\[
R(z)=I+z_1\lambda(a)+z_2\lambda(t).
\]
They showed that
\[
\log\det R(z)
=\frac1{4\pi}\int_0^{2\pi}
\log|1-z_1^2-z_2^2-2z_1z_2\cos s|\,ds
\]
and related the traced Maurer--Cartan form $R^{-1}dR$ to the first homology group of the projective resolvent set \cite[Sec.~4]{GY}. Further determinant formulas and their relation to the geometry of the resolvent set are discussed in \cite[Sec.~5.3]{YangBook}.

Let $(\mathcal M,\tau)$ be a finite von Neumann algebra. Let $M_d(\mathcal M)$ be the algebra of $d\times d$ matrices over $\mathcal M$, equipped with the unnormalized trace $\Tr_d\otimes\tau$.
\begin{introdefinition}\label{d1}
For $\mathbf{T}\in M_d(\mathcal M)$, its \emph{Fuglede--Kadison determinant} is
\[
\det (\mathbf{T})
=\exp\left(\int_{[0,\infty)}\log s\,d\mu_{|\mathbf{T}|}(s)\right),
\ \ |\mathbf{T}|=(\mathbf{T}^*\mathbf{T})^{1/2}.
\]
The von Neumann algebra and its trace are understood from the context. With a bit abuse of notation, The notation $\det \mathbf{T}$ is used for both the FK-determinant in a finite von Neumann algebra and the ordinary determinant of a complex matrix, the meaning will be clear from the context. The meaning should also be clear from the context. The measure $\mu_{|\mathbf{T}|}$ is defined via the spectral projections of $|\mathbf{T}|$ by
\[
\mu_{|\mathbf{T}|}(E)=(\Tr_d\otimes\tau)\bigl(\mathbf \chi_E(|\mathbf{T}|)\bigr)
\]
for every Borel set $E\subset[0,\infty)$, where $\chi_E$ stands for the characteristic function of set $E$. For $d=1$, this agrees with the FK-determinant on $\mathcal M$ \cite{FK}.
\end{introdefinition}
For an invertible differentiable function $A: \C^m\to {\mathcal M}$, the Jacobi formula implies
\[
d\log\det A(z)
=\operatorname{Re}\tau\bigl(A(z)^{-1}dA(z)\bigr),
\]
which relates the determinant to the traced Maurer-Cartan form \cite[Sec.~5.3]{YangBook}. For every discrete group $\mathscr{G}$, the group von Neumann algebra $L(\mathscr{G})$ has a normal trace $\tau_\mathscr{G}$ defined by 
\[\tau_\mathscr{G}(a):=\langle a\delta_{\bf 1}, \delta_{\bf 1}\rangle,\] and thus the FK-determinant is also well-defined on invertible elements $a\in L(\mathscr{G})$.

Atiyah's $L^2$-index theorem introduced von Neumann dimension into index theory for coverings \cite{Atiyah}. FK-determinants arise in the theory of $L^2$-torsion and determinant class developed by L\"uck and Schick \cite{Lueck,SchickDet}; their relation to $K$-theory and Whitehead torsion is surveyed in \cite{deLaHarpe}. In algebraic dynamics, Lind--Schmidt--Ward related Mahler measure to entropy for $\mathbb{Z}^m$-actions \cite{LSW}. Deninger connected FK-determinants with entropy and Mahler measure, and expressed the determinants of crossed-product operators in terms of Lyapunov exponents \cite{DeningerEntropy,DeningerMahler,DeningerLyap}.

Li--Thom proved that, for a countable discrete amenable group $\Gamma$ and $f\in\mathbb{Z}\Gamma$ with $\ker\lambda(f)=\{0\}$, the principal algebraic action on
\[
X_f=\widehat{\mathbb{Z}\Gamma/\mathbb{Z}\Gamma f}
\]
satisfies
\[
h(\Gamma\curvearrowright X_f)
=\log\det\lambda(f),
\]
confirming Deninger's conjecture \cite[Theorem~1.2]{LiThom}. Their entropy--$L^2$-torsion theorem also implies the vanishing of $L^2$-torsion for nontrivial amenable groups with a finite classifying space \cite[Theorems~1.1 and~1.3]{LiThom}. A key analytic ingredient is a Szeg\H{o}-type formula that approximates FK-determinants of positive operators by normalized determinants of finite compressions along F{\o}lner sets \cite[Theorem~1.4]{LiThom}. Hayes later obtained entropy--determinant formulas for algebraic actions of sofic groups \cite{Hayes}.

The main analytic difficulty comes from small singular values. Injectivity rules out an atom at zero, whereas positivity of the FK-determinant additionally requires
\[
\int_{(0,1)}|\log s|\,d\mu_{|\mathbf{T}|}(s)<\infty.
\]
For amenable groups, Li--Thom proved positivity for injective matrices with algebraic coefficients and asked the following question \cite[Lemma~3.12 and Question~3.13]{LiThom}.

\medskip

\noindent\textbf{Injectivity and determinant positivity question}. Do we have $\det(\lambda(\mathbf{A}))>0$ for all $\mathbf{A}\in M_d(\C\Gamma)$ with $\ker\lambda(\mathbf{A})=\{0\}$?

\medskip

For any countable group $\mathscr{G}$, an element $\mathfrak e\in \mathbb{C}[\mathscr{G}]$ is called a projection if $\mathfrak e=\mathfrak e^*=\mathfrak e^2$. Let $K$ be a discrete group (not necessarily abelian) whose group algebra $\C[K]$ contains a nontrivial projection $\mathfrak e$. We shall study the FK-determinant of the pencil
\[
\mathfrak f_{\kappa,\beta,\gamma}
=\kappa\one+\beta t^{-1}\mathfrak e+\gamma \mathfrak et
\]
in the left regular representation of $K\wrp\Z$. The orthogonal decomposition of Dicks--Schick reduces this operator to weighted finite Jacobi matrices \cite{DicksSchick}. A continued-fraction construction yields injective scalar elements with vanishing FK-determinant. On the other hand, let $G_{tf}\rtimes_\sigma\mathbb Z$, where $G_{tf}$ is any countable torsion-free
abelian group and $\sigma\in\operatorname{Aut}(G_{tf})$, the Jensen--Mahler estimate together with Ore localization yield determinant positivity for matrices of arbitrary size. This paper establishes the following facts.

\begin{theorem}\label{t3}
\begin{enumerate}[label=\textup{(\roman*)}]
\item Assume $W:=\tau_K(\mathfrak e)^{-1}>1$. There exists $\alpha_*\in\R\setminus\Q$, depending only on $W$, such that for every $\beta,\gamma\in\C^\times$, every $\upsilon$ with $\upsilon^2=\beta\gamma$, and $\kappa=2\upsilon\cos(\pi\alpha_*)$, we have
\[
\ker\lambda(\mathfrak f_{\kappa,\beta,\gamma})=\{0\},\ \
\det\lambda(\mathfrak f_{\kappa,\beta,\gamma})=0.
\]

\item Let $\Gamma=G_{tf}\rtimes_\sigma\mathbb Z$, where $G_{tf}$ is any countable torsion-free
abelian group and $\sigma\in\operatorname{Aut}(G_{tf})$. For every $d\ge1$ and every $\mathbf{A}\in M_d(\C\Gamma)$, injectivity of $\lambda(\mathbf{A})$ implies $\det\lambda(\mathbf{A})>0$.
\end{enumerate}
\end{theorem}

If $K$ is countable amenable and contains a nontrivial finite subgroup, then the averaging projection onto that subgroup satisfies the hypothesis of part~\textup{(i)}. Since $K\wr\mathbb{Z}$ is amenable, part~\textup{(i)} provides a negative answer to the Li--Thom question for $\Gamma=K\wr\mathbb{Z}$. In particular, this applies to every nontrivial finite group $\mathscr{G}$. Part~\textup{(ii)} gives an affirmative answer for $G_{tf}\rtimes_\sigma\mathbb Z$ in every matrix size. These results do not yield a general dichotomy between torsion and torsion-free groups.

The paper is organized as follows. Section~\ref{s2} introduces the group-theoretic notation, the Koopman and left regular representations, and the recurrence formulas for the finite Jacobi matrices used later. Section~\ref{s3} studies the Schur recursion, proves weak equivalence of the Koopman and regular representations, and determines the parts of the projective spectrum lying in the Julia and Fatou sets. Section~\ref{s4} applies the Dicks--Schick decomposition to the relevant wreath-product operators, computes their FK-determinants, and proves Theorem~\ref{t3}.

\section{Preliminaries}\label{s2}

Throughout Sections~\ref{s2}--\ref{s3}, $G$ denotes a nontrivial finite abelian group of order $r\ge2$, and
$\Gamma_G=G\wrp\Z$.
We use the notation in \eqref{e3}--\eqref{e7}.

\subsection{The wreath product and the group-ring pencil}

Let
$
\mathcal{G}=\bigoplus\limits_{j\in\Z}G_j,
$
where $G_j$ is a copy of $G$ in the $j$th coordinate.  Let $x=(x_j)_{j\in\Z}\in \mathcal{G}$, and $(\vartheta x)_j=x_{j-1}$.  Then
\[
\Gamma_G=\mathcal{G}\rtimes_\vartheta\Z,\ \
(x,n)(y,m)=\bigl(x\vartheta^n(y),n+m\bigr).
\]
We identify $G$ with $G_0$.  The shift generator $t$ introduced before \eqref{e3} is $t=(\one,1)$, hence
\[
tbt^{-1}=\vartheta(b),\ \ b\in \mathcal{G}.
\]

The transformations $s_j$ are given by \eqref{e3}.  They generate $\Gamma_G$.  Indeed, $s_1=t$ and $s_js_1^{-1}=g_j^{-1}$, so they generate $G_0$ and $t$.  Their conjugates generate every $G_j$, and hence $\mathcal{G}\rtimes\langle t\rangle$.  
Since the group-ring involution is $g^*=g^{-1}$, hence
\[
s_j^*=t^{-1}g_j.
\]
With $e$ and $C_G$ as in \eqref{e5}, define
\begin{equation}\label{e8}
X=rY=rt^{-1}e.
\end{equation}
Then
\begin{equation}\label{e9}
\sum_{j=1}^r s_j^*=X,\ \
\sum_{j=1}^r s_j=X^*.
\end{equation}
For $z=(z_0,z_1,z_2,z_3)\in\C^4$, let
\begin{equation}\label{e10}
D_G(z)=z_0\one+z_1C_G+z_2X+z_3X^*\in\C[\Gamma_G].
\end{equation}
For a unitary representation $\rho:\Gamma_G\to\mathcal U(\mathscr H_\rho)$, where $\mathscr H_\rho$ is the Hilbert space on which $\rho$ acts. 
Let $D_{\rho,G}(z)=\rho(D_G(z))$.  Since $D_G(az)=aD_G(z)$ for $a\ne0$, its invertibility depends only on $[z]\in\Pj^3$.

We shall use the following elementary property of $Y=t^{-1}e$. 
\begin{lemma}\label{l1}
For every unitary representation $\rho$ of $\Gamma_G$ and every $k\ge1$, $\rho(Y)^k$ is a partial isometry.  More precisely, if
\[
e_j=t^{-j}et^j,\ \
E_k=e_{-(k-1)}\cdots e_{-1}e_0,
\]
then
\[
Y^k=t^{-k}E_k,
\]
and $E_k$ is a projection.
\end{lemma}

\begin{proof}
For each $j$, the element
\[
e_j=\frac1r\sum\limits_{g\in G_{-j}}g
\]
is the average over the finite subgroup $G_{-j}$ and is an orthogonal projection in every unitary representation.  If $i\ne j$, then $G_{-i}$ and $G_{-j}$ commute in the base group $\mathcal{G}=\bigoplus_{m\in\Z}G_m$.  Hence
\[
e_ie_j
=\frac1{r^2}\sum\limits_{g\in G_{-i}}\sum\limits_{h\in G_{-j}}gh
=\frac1{r^2}\sum\limits_{g\in G_{-i}}\sum\limits_{h\in G_{-j}}hg
=e_je_i.
\]
Thus the $e_j$ are pairwise commuting projections, so $E_k$ is a projection.  Since $et^{-1}=t^{-1}e_{-1}$, and by induction, we have $Y^k=t^{-k}E_k$.  The product of a unitary and a projection is a partial isometry.  Applying $\rho$ yields the result.
\end{proof}

\subsection{The Koopman and left regular representations}

Recall from \eqref{e4} and \eqref{e6} that
$\mathscr H_n=\ell^2(G^n)$, $X_n=\sum_{j=1}^r\pi_n(s_j)^*$, and $S_n=\pi_n(C_G)$.  From \eqref{e8},
\[
Y_n=\frac1rX_n.
\]
By \eqref{e8}--\eqref{e9}, then
$Y_n=\pi_n(Y)$, $S_n=\pi_n(C_G)$, $\pi_n(D_G(z))=D_n(z)$,
where $D_n(z)$ is the pencil defined in \eqref{e7}.  We use the identity
\begin{equation}\label{e11}
X_n^*X_n=r(I+S_n),\ \
Y_n^*Y_n=\frac1r(I+S_n),
\end{equation}
which follows from \cite[Lemma~5.1]{KSS}.

Since $G^{n+1}=G\times G^n$, we have
\[
\mathscr H_{n+1}=\ell^2(G^{n+1})
\cong \ell^2(G)\otimes \mathscr H_n
\cong \C^r\otimes \mathscr H_n.
\]
Let $J_r$ be the $r\times r$ all-ones matrix.  The operator
\begin{equation}\label{e12}
P_n=\frac1rJ_r\otimes I_{\mathscr H_n}:\mathscr H_{n+1}\longrightarrow \mathscr H_{n+1}
\end{equation}
is the orthogonal projection onto
\[
P_n\mathscr H_{n+1}
=\C(1,\ldots,1)^T\otimes \mathscr H_n.
\]
Let $Q_n=I_{\mathscr H_{n+1}}-P_n$.  The map
\begin{equation}\label{e13}
V_n:\mathscr H_n\longrightarrow P_n\mathscr H_{n+1},\ \
V_n\xi=r^{-1/2}(1,\ldots,1)^T\otimes\xi,
\end{equation}
is unitary.  Thus $V_n$ maps $\mathscr H_n$ unitarily onto $\C(1,\ldots,1)^T\otimes \mathscr H_n$.

Equip $G^{\N}$ with the product probability measure $\nu$. The action on $(G^{\N},\nu)$ yields the Koopman representation \eqref{e1}. In the interval model described in the Introduction, the cylinder functions depending on the first $n$ coordinates form the reducing subspace $\mathscr H_n$ defined in \eqref{e4}. The subspaces $\mathscr H_n$ form an increasing sequence whose union is dense in $\mathscr H=L^2(G^{\N},\nu)$. Let
\[
D_{\pi,G}(z)=\pi(D_G(z)):\mathscr H\to \mathscr H.
\]

\begin{proposition}\label{p1}
One has
\[
D_{\pi,G}(z)\text{ invertible}
\iff
\left\{
\begin{array}{l}
\forall n\ge0,\ \ D_n(z)\text{ invertible},\\[1mm]
\displaystyle\sup_{n\ge0}\norm{D_n(z)^{-1}}<\infty.
\end{array}
\right.
\]
\end{proposition}

\begin{proof}
If $D_{\pi,G}(z)$ is invertible, its restriction to each reducing subspace $\mathscr H_n$ is $D_n(z)$, and
$\norm{D_n(z)^{-1}}\le \norm{D_{\pi,G}(z)^{-1}}$.

Conversely, assume
\[
\forall n\ge0, \ \ D_n(z)\ \ \text{invertible}, \ \ \sup\limits_{n\ge0}\|D_n(z)^{-1}\|<\infty.
\]
Since $\mathscr H_n\subset \mathscr H_{n+1}$ is a reducing subspace for $D_{\pi,G}(z)$, we  have
\[
D_{n+1}(z)^{-1}\big|_{\mathscr H_n}=D_n(z)^{-1},\ \ \forall n\ge0.
\]
Hence $\{D_n(z)^{-1}\}_{n\ge0}$ defines a bounded operator on $\bigcup\limits_{n\ge0}\mathscr H_n$.  Since this union is dense in $\mathscr H$, the operator extends uniquely to a bounded inverse of $D_{\pi,G}(z)$.
\end{proof}

Let $\lambda$ denote the left regular representation of $\Gamma_G$, obtained from \eqref{e2} by taking $K=\Gamma_G$.  Let
\[
\Nalg\Gamma_G=\lambda(\Gamma_G)'', \ \
\tau_{\Gamma_G}(T)=\langle T\delta_{1_{\Gamma_G}},\delta_{1_{\Gamma_G}}\rangle, \ \
T \in \Nalg\Gamma_G,
\]
and let
\[
D_{\lambda,G}(z)=\lambda(D_G(z)):\ell^2(\Gamma_G)\to\ell^2(\Gamma_G).
\]
We denote the corresponding noninvertibility sets in $\Pj^3$ by
$p(D_{\pi,G})$ and $p(D_{\lambda,G})$.

\subsection{Finite Jacobi matrices and recurrence polynomials}

For $z=(z_0,z_1,z_2,z_3)\in\C^4$, let
\[
\theta=z_0+(r-1)z_1,\ \
\delta=z_0-z_1,\ \
u=z_2,\ \
v=z_3,\ \
c=r^2uv.
\]
We use these quantities throughout Sections~\ref{s3}--\ref{s4}.

For $a,b\in\C$, define
\begin{equation}\label{e14}
\mathsf T_{-1}(a,b)=0,\ \
\mathsf T_0(a,b)=1,\ \
\mathsf T_{m+1}(a,b)=a\mathsf T_m(a,b)-b\mathsf T_{m-1}(a,b).
\end{equation}
Equivalently,
\[
\mathsf T_m(a,b)
=\sum\limits_{j=0}^{\lfloor m/2\rfloor}(-b)^j
\binom{m-j}{j}a^{m-2j}.
\]

\noindent\textbf{The Jacobi matrices $\mathcal{A}_m(z)$.}
For $m\ge1$, let $J_m$ be the nilpotent backward shift on $\C^m$ with respect to the coordinate orthonormal basis $e_0,\ldots,e_{m-1}$:
\[
J_me_0=0,\ \
J_me_j=e_{j-1}\ \
(1\le j\le m-1).
\]
In the terminology of Halmos--Wallen, $J_m$ is a truncated shift of index $m$ \cite[Theorem~1]{HalmosWallen}.  We use this terminology in what follows.
Let
\begin{equation}\label{e15}
\mathcal{A}_m(z)=
\begin{pmatrix}
\delta&ru&&&\\
rv&\theta&ru&&\\
&rv&\theta&\ddots&\\
&&\ddots&\ddots&ru\\
&&&rv&\theta
\end{pmatrix}.
\end{equation}
The matrices \eqref{e15} appear in both representations considered below. In the Koopman representation, they arise from the decomposition of $Y_n$ into truncated shifts in Theorem~\ref{t5}. In the left regular representation, Dicks--Schick construct mutually orthogonal subspaces spanned by the vectors $t^j q_n$ and represent the operator on each such subspace by a tridiagonal matrix, which can be found in \cite[(3.3)--(3.5), Lemmas~3.4 and~3.6, and Remark~3.7]{DicksSchick}. Since $C_G=rY^*Y-I$, the restriction of the pencil is of the form \eqref{e15}. Theorem~\ref{t9} recalls the orthogonality, completeness, and von Neumann weights needed in Section~\ref{s4}.

Define
\begin{equation}\label{e16}
\mathsf G_{-1}=1,\ \
\mathsf G_n=\delta\mathsf T_n(\theta,c)-c\mathsf T_{n-1}(\theta,c).
\end{equation}
Then
\begin{equation}\label{e17}
\mathsf G_0=\delta,\ \
\mathsf G_{n+1}=\theta\mathsf G_n-c\mathsf G_{n-1}.
\end{equation}
Expansion along the first row of \eqref{e15} implies that
\begin{equation}\label{e18}
\det \mathcal{A}_m(z)=\mathsf G_{m-1}(z).
\end{equation}

Let $x=(x_1,\ldots,x_m)^T\in\C^m$. For $2\le j\le m-1$, the $j$th interior row equation of $\mathcal{A}_m(z)x=0$ is
$rvx_{j-1}+\theta x_j+ru x_{j+1}=0$.
For a solution of the form $x_j=\xi^{j-1}$, we obtain
$ru\xi^2+\theta\xi+rv=0$.
With $\varrho=-ru\xi$, then
\begin{equation}\label{e19}
\varrho^2-\theta\varrho+r^2uv=0.
\end{equation}
If $\varrho_1,\varrho_2$ are its roots, one has
\begin{equation}\label{e20}
\mathsf T_m(\theta,c)=
\begin{cases}
\displaystyle\frac{\varrho_1^{m+1}-\varrho_2^{m+1}}{\varrho_1-\varrho_2},
&\varrho_1\ne\varrho_2,\\[3mm]
(m+1)\varrho_1^m,&\varrho_1=\varrho_2.
\end{cases}
\end{equation}
These are the standard formulas for the characteristic roots of a constant-coefficient second-order difference equation, we refer the reader to \cite[Secs.~1.1--1.2]{Teschl} for details. For the construction of the Green function from solutions decaying in both directions, see \cite[Sec.~1.2, especially Eq.~(1.99)]{Teschl}.

\noindent\textbf{The Jacobi matrices $\mathcal{B}_m$.}
For $m\ge1$ and $a,\alpha,\beta\in\C$, let
\[
\mathcal{B}_m(a;\alpha,\beta)=
\begin{pmatrix}
a&\alpha&&&\\
\beta&a&\alpha&&\\
&\beta&a&\ddots&\\
&&\ddots&\ddots&\alpha\\
&&&\beta&a
\end{pmatrix}.
\]
Expansion along the last row shows that
\begin{equation}\label{e21}
\det \mathcal{B}_m(a;\alpha,\beta)=\mathsf T_m(a,\alpha\beta).
\end{equation}
If $b\ne0$, then
\begin{equation}\label{e22}
\mathsf T_m(a,b)
=(\sqrt b)^m U_m\!\left(\frac{a}{2\sqrt b}\right),
\end{equation}
where $U_m$ is the Chebyshev polynomial of the second kind.  The expression in \eqref{e22} is independent of the choice of $\sqrt b$.  Thus
\[
U_{-1}=0,\qquad U_0=1,\qquad U_{m+1}(x)=2xU_m(x)-U_{m-1}(x),
\]
and
\[
U_m(\cos t)=\frac{\sin((m+1)t)}{\sin t},
\qquad
U_m(\cosh s)=\frac{\sinh((m+1)s)}{\sinh s},
\]
with the usual limiting interpretation when the denominator vanishes.  See \cite[Sec.~18.5]{DLMF}.

\section{Self-similar spectral dynamics and weak equivalence}\label{s3}

\subsection{Schur reduction and the projective spectra}\label{ss1}

Recall the operators $P_n,Q_n$ and the unitary $V_n$ from \eqref{e12}--\eqref{e13}.  In particular,
\begin{equation}\label{e23}
\mathscr H_{n+1}=P_n\mathscr H_{n+1}\oplus Q_n\mathscr H_{n+1},
\end{equation}
and $P_n\mathscr H_{n+1}\cong \mathscr H_n$ through $V_n$.  
Define
\[
\mathcal U_n=\diag\bigl(\pi_n(s_1),\ldots,\pi_n(s_r)\bigr)
\in\mathcal B(\C^r\otimes \mathscr H_n).
\]
By \cite[Sec.~3, especially (3.6)--(3.9), and Sec.~5.1]{KSS},
\begin{equation}\label{e24}
\frac1r\sum_{i=1}^r\pi_{n+1}(s_i)=P_n\mathcal U_n,\ \
Y_{n+1}=\mathcal U_n^*P_n.
\end{equation}
Also we have $\mathcal U_n=\pi_{n+1}(t)$ and $\pi_{n+1}(e)=P_n$.

Kambites--Silva--Steinberg prove
$\sum\limits_{i\ne j}\pi_n(s_i)\pi_n(s_j)^*=rS_n$,
\cite[Lemma~5.1]{KSS}. Hence
\[
X_n^*X_n=r(I+S_n),\ \
Y_n^*Y_n=\frac1r(I+S_n),\ \
S_{n+1}=rP_n-I.
\]

The variables introduced in Subsection~2.3 determine $z$ uniquely:
\[
z_0=\frac{\theta+(r-1)\delta}{r},\ \
z_1=\frac{\theta-\delta}{r},\ \
z_2=u,\ \
z_3=v.
\]
For $\delta\ne0$, in $\C^4$, define
\begin{equation}\label{e25}
F(z)=\left(
\theta-\frac{r(r-1)uv}{\delta},\,
\frac{ruv}{\delta},\,
u,\,
v
\right).
\end{equation}
Since
\[
\theta(F(z))=\theta,\ \
u(F(z))=u,\ \
v(F(z))=v,\ \
\delta(F(z))=\theta-\frac{r^2uv}{\delta},
\]
the same map can be written as
\begin{equation}\label{e26}
F(\theta,\delta,u,v)
=\left(\theta,\theta-\frac{r^2uv}{\delta},u,v\right).
\end{equation}

Let $z^{(0)}=z$ and $\delta_0=\delta$.  Recursively, whenever $\delta_j\ne0$, let
\begin{equation}\label{e27}
z^{(j+1)}=F(z^{(j)}),\ \
\delta_{j+1}=(z^{(j+1)})_0-(z^{(j+1)})_1.
\end{equation}
Thus $F^k(z)$ is well defined when $\delta_j\ne0$, $0\le j<k$.
For $n=0$, we have $\mathscr H_0=\C$ and
\begin{equation}\label{e28}
D_0(z)=\theta+r(u+v).
\end{equation}

For
\[
M=\begin{pmatrix}A&B\\ C&D\end{pmatrix},
\  \ with \ D \ invertible,
\]
the Schur complement is defined by
\[
S_D=A-BD^{-1}C.
\]
Then $M$ is invertible if and only if $S_D$ is invertible, and $\det M=(\det D)(\det S_D)$.

\begin{proposition}\label{p2}
Assume $\delta\ne0$,  with respect to the decomposition \eqref{e23}, the Schur complement of $\delta Q_n$ in $D_{n+1}(z)$ is unitarily equivalent to $D_n(F(z))$.  Consequently,
\begin{equation}\label{e29}
\det D_{n+1}(z)
=\delta^{(r-1)r^n}\det D_n(F(z)),
\end{equation}
and $D_{n+1}(z)$ is invertible if and only if $D_n(F(z))$ is invertible.  If $\delta_0,\ldots,\delta_n$ are nonzero, then
\begin{equation}\label{e30}
\det D_{n+1}(z)
=
\left(\prod_{j=0}^{n}\delta_j^{(r-1)r^{\,n-j}}\right)
D_0\bigl(z^{(n+1)}\bigr),
\end{equation}
where $D_0$ is the scalar pencil in \eqref{e28}.
\end{proposition}

\begin{proof}
From
\[
S_{n+1}=rP_n-I,\ \
X_{n+1}=r\mathcal U_n^*P_n,\ \
X_{n+1}^*=rP_n\mathcal U_n,
\]
we obtain
\[
D_{n+1}(z)=\delta I+rz_1P_n+ru\mathcal U_n^*P_n+rvP_n\mathcal U_n.
\]
With respect to \eqref{e23}, this is
\begin{equation}\label{e31}
	D_{n+1}(z)=
	\begin{pmatrix}
		\theta P_n+ruP_n\mathcal U_n^*P_n+rvP_n\mathcal U_nP_n
		&rvP_n\mathcal U_nQ_n\\[1mm]
		ruQ_n\mathcal U_n^*P_n
		&\delta Q_n
	\end{pmatrix}.
\end{equation}
Since $\delta\ne0$, $\delta Q_n$ is invertible on $Q_n\mathscr H_{n+1}$.  Hence
\[
\begin{pmatrix}I&-rvP_n\mathcal U_nQ_n(\delta Q_n)^{-1}\\0&I\end{pmatrix}
D_{n+1}(z)
\begin{pmatrix}I&0\\-ru(\delta Q_n)^{-1}Q_n\mathcal U_n^*P_n&I\end{pmatrix}
=
\begin{pmatrix}\mathcal T_n&0\\0&\delta Q_n\end{pmatrix},
\]
where
\[
\mathcal T_n=\theta P_n+ruP_n\mathcal U_n^*P_n+rvP_n\mathcal U_nP_n-r^2uvP_n\mathcal U_nQ_n(\delta Q_n)^{-1}Q_n\mathcal U_n^*P_n.
\]
Substituting the blocks from \eqref{e31} implies that
\[
\mathcal T_n
=\theta P_n+ruP_n\mathcal U_n^*P_n+rvP_n\mathcal U_nP_n
-\frac{r^2uv}{\delta}P_n\mathcal U_nQ_n\mathcal U_n^*P_n.
\]
By \eqref{e13}, then
$P_nV_n=V_n$ and $V_n^*P_n=V_n^*$.  Since
\[
\mathcal U_n^*=\diag\bigl(\pi_n(s_1)^*,\ldots,\pi_n(s_r)^*\bigr),
\]
for $\xi\in \mathscr H_n$, we have
\[
V_n^*P_n\mathcal U_n^*P_nV_n\xi
=V_n^*\mathcal U_n^*V_n\xi
=\frac1r\sum\limits_{j=1}^r\pi_n(s_j)^*\xi
=Y_n\xi.
\]
Thus
\[
V_n^*P_n\mathcal U_n^*P_nV_n=Y_n,\ \ 
V_n^*P_n\mathcal U_nP_nV_n=Y_n^*.
\]
From $Q_n=I-P_n$, the unitarity of $\mathcal U_n$, and \eqref{e11}, we have
\[
V_n^*P_n\mathcal U_nQ_n\mathcal U_n^*P_nV_n
=I-Y_n^*Y_n
=\frac{r-1}{r}I-\frac1rS_n.
\]
Therefore,
\[
V_n^*\mathcal T_nV_n
=\left(\theta-\frac{r(r-1)uv}{\delta}\right)I
+\frac{ruv}{\delta}S_n+uX_n+vX_n^*
=D_n(F(z)).
\]
Since
$\dim Q_n\mathscr H_{n+1}=(r-1)r^n$, we obtain \eqref{e29}.  Iterating \eqref{e29} yields \eqref{e30}.
\end{proof}

For every $\lambda\in\C^\times$, $F(\lambda z)=\lambda F(z)$.  Therefore $F$ induces the rational map
\begin{equation}\label{e32}
\widehat F:\Pj^3\dashrightarrow\Pj^3,
\qquad
\widehat F([\theta:\delta:u:v])
=[\theta\delta:\theta\delta-r^2uv:u\delta:v\delta].
\end{equation}
From \eqref{e26}, we have
\[
\delta_{n+1}=\theta-\frac{c}{\delta_n}.
\]
Combining this recursion with \eqref{e17} and \eqref{e18}, whenever
$\delta_j\ne0$ for $0\le j<n$, we have
\begin{equation}\label{e33}
\delta_n
=\frac{\mathsf G_n}{\mathsf G_{n-1}}
=\frac{\det \mathcal{A}_{n+1}(z)}{\det \mathcal{A}_n(z)}.
\end{equation}
Consequently,
\[
\det \mathcal{A}_m(z)=\prod_{j=0}^{m-1}\delta_j
\]
whenever $\delta_j\ne0$ for $0\le j<m$.  In particular, if
$\delta_j\ne0$ for $0\le j<n$, then
\[
\mathsf G_n(z)=0
\Longleftrightarrow
\delta_n=0
\Longleftrightarrow
\det \mathcal{A}_{n+1}(z)=0.
\]

\begin{remark}
Proposition~\ref{p2} shows that two different issues enter the projective spectrum. A matrix $D_n(z)$ becomes singular when $\delta_j(z)=0$ for some $j\ge0$ or when an iterate reaches the hyperplane $L_G$ corresponding to $n=0$. These singularities do not describe the whole projective spectrum: since the matrices $D_n(z)$ are not normal in general, every $D_n(z)$ may be invertible while $\sup\limits_{n\geq 1}\|D_n(z)^{-1}\|=\infty$. For a normal matrix, the smallest singular value equals the distance from 0 to its spectrum, whereas for a nonnormal matrix small singular values may occur far from the eigenvalues, this is the basic pseudospectral phenomenon discussed in \cite[Secs.~2--5]{TrefethenEmbree}. We therefore consider both $\delta_j=0$ and the possible divergence of $\|D_n(z)^{-1}\|$ as $n\to\infty$.
\end{remark}

Repeated use of \eqref{e29} reduces the determinant to $D_0$.  We first determine its zero set.

\noindent\textbf{The hyperplane $L_G$} \ \
By \eqref{e28},
$
D_0(z)=\theta+r(u+v).
$
Let $\mathbf 1_n$ be the constant vector in $\mathscr H_n$.  Then
\[
S_n\mathbf 1_n=(r-1)\mathbf 1_n,\ \
D_n(z)\mathbf 1_n=D_0(z)\mathbf 1_n.
\]
Thus $\C\mathbf 1_n$ and $\mathscr H_n^0=\mathbf 1_n^\perp$ reduce $D_n(z)$, and let
\begin{equation}\label{e34}
\begin{aligned}
L_G
&=\{[\theta:\delta:u:v]\in\Pj^3\mid \theta+r(u+v)=0\}\\
&=\left\{[z_0:z_1:z_2:z_3]\in\Pj^3\mid
z_0+(r-1)z_1+r(z_2+z_3)=0\right\}.
\end{aligned}
\end{equation}
By \eqref{e26},
\begin{equation}\label{e35}
D_0(F^k(z))=D_0(z)\ \
\bigl(k\ge0,\ \delta_j(z)\ne0\ (0\le j<k)\bigr).
\end{equation}
Equivalently,
\begin{equation}\label{e36}
\left\{[z]\middle|
\delta_j(z)\ne0\ (0\le j<k),\
\widehat F^k([z])\in L_G
\right\}
=
\left\{[z]\middle|
\delta_j(z)\ne0\ (0\le j<k),\
[z]\in L_G
\right\}.
\end{equation}
Moreover,
\begin{equation}\label{e37}
D_{\pi,G}(z)\mathbf 1=D_0(z)\mathbf 1,\ \
L_G\subset p(D_{\pi,G}).
\end{equation}

By \eqref{e26}, $\theta,u,v$ remain unchanged, while $\delta$ is replaced by $\theta-r^2uv/\delta$.  Hence
\[
\Pi:\Pj^3\dashrightarrow\Pj^2,\ \
\Pi([\theta:\delta:u:v])=[\theta:u:v],\ \
\Pi\circ\widehat F=\Pi,
\]
where the rational projection $\Pi$ is undefined at $[0:1:0:0]$, and $\Pi^{-1}(\eta)$ is an affine line.
For $\eta=[\theta:u:v]\in\Pj^2$, let
\[
L_\eta=\Pi^{-1}(\eta).
\]
Away from the indeterminacy points, the map $\widehat F$ preserves the value of $L_\eta$
and the induced map on $L_\eta$ is given by
\[
M_\eta(\delta)=\theta-\frac{c}{\delta},\ \
c=r^2uv.
\]

For $uv\ne0$, $M_\eta$ extends to a Möbius transformation on $\mathbb P^1$. Let $\mathcal F_\eta$ and $\mathcal J_\eta$ denote its Fatou and Julia sets.
Let $I(\widehat F^n)$ be the indeterminacy set of $\widehat F^n$ and let
\[
\mathcal D_{n+1}
=
\left\{
[z]\in\mathcal D_n:
\widehat F^{n}([z])\notin I(\widehat F)
\right\},\ \
\mathcal D_0=\mathbb P^3,\ \
\Omega_\infty=\bigcap_{n\ge1}\mathcal D_n.
\]
For every $n\ge0$, the map $\widehat F^n$ is holomorphic on $\Omega_\infty$.

A family of holomorphic maps on $U$ is normal if every sequence in the family has a subsequence converging uniformly on compact subsets of $U$.  Since $\Pj^3$ is compact, normality is equivalent to local equicontinuity in local holomorphic coordinates.

\begin{definition}
\[
\mathcal F(\widehat F)
=
\left\{
[z]\in\Omega_\infty
\middle|
\begin{array}{c}
\exists\text{ an open neighborhood }U\ni [z],\ U\subset\Omega_\infty,\\
\{\widehat F^n([z])|_U\}_{n\ge0}\text{ is normal}
\end{array}
\right\},\ \
\mathcal J(\widehat F)=\Pj^3\setminus\mathcal F(\widehat F).
\]
\end{definition}

\begin{remark}\label{a1}
Equip $\mathbb P^3$ and $\mathbb P^1$ with their usual topologies. In projective space, convergence is equivalent to convergence of affine coordinates in any chart containing the limit.
For the affine lift $F$, projective convergence also implies convergence of suitably normalized homogeneous representatives. Indeed, let
$z\in\mathbb C^4\setminus\{0\}$ be such that all iterates $F^k(z)$ are well defined. Let
\[
\theta=z_0+(r-1)z_1,\ \
\delta_k=F^k(z)_0-F^k(z)_1,\  \
u=z_2,\ \ v=z_3.
\]
Since $F$ preserves $\theta,u,v$, we have
\[
F^k(z)=
\left(
\frac{\theta+(r-1)\delta_k}{r},
\frac{\theta-\delta_k}{r},u,v
\right).
\]
Suppose that
$\lim\limits_{k\to\infty}\widehat F^{\,k}([z])=\xi$, $\xi\in E_{1,\mathrm{att}}$.
Then $\xi$ has a unique representative
$\zeta\in\mathbb C^4\setminus\{0\}$ such that
\[
\zeta_0+(r-1)\zeta_1=\theta,\ \
\zeta_2=u,\ \ \zeta_3=v.
\]
Let $\delta_\infty=\zeta_0-\zeta_1$, projective convergence in any affine chart in which at least one of $\theta,u,v$ is nonzero yields
$\delta_k\to\delta_\infty$. Consequently,
\[
\lim\limits_{k\to\infty}F^{\,k}(z)=\zeta
=
\left(
\frac{\theta+(r-1)\delta_\infty}{r},
\frac{\theta-\delta_\infty}{r},u,v
\right)
\ \ \text{in }\mathbb C^4.
\]
For estimates involving homogeneous vectors,
$\|\cdot\|_2$ denotes the Euclidean norm on $\mathbb C^4$.
\end{remark}

Let
\[
\Delta_\eta=\theta^2-4r^2uv,\ \
\Omega_{\mathrm{reg}}
=
\left\{
[\theta:\delta:u:v]\in\Pj^3
\middle|
uv\ne0, \Delta_\eta\ne0
\right\}.
\]
For $[z]\in\Omega_{\mathrm{reg}}$, select a representative
$z=(\theta,\delta,u,v)$.
Let $\varrho_1(z),\varrho_2(z)$ be two roots of
\begin{equation}\label{e38}
\varrho^2-\theta\varrho+r^2uv=0.
\end{equation}
Let
\begin{equation}\label{e39}
\psi=\frac{\varrho_2}{\varrho_1},\ \
w=\frac{\delta-\varrho_1}{\delta-\varrho_2}.
\end{equation}
Since $\varrho_1+\varrho_2=\theta$ and $\varrho_1\varrho_2=r^2uv$,
\begin{equation}\label{e40}
\frac{M_\eta(\delta)-\varrho_1}{M_\eta(\delta)-\varrho_2}
=\psi\frac{\delta-\varrho_1}{\delta-\varrho_2},\ \
\widehat F^n(\eta,w)=(\eta,\psi^nw).
\end{equation}
Moreover,
\[
M_\eta'(\varrho_1)=\frac{r^2uv}{\varrho_1^2}=\frac{\varrho_2}{\varrho_1}=\psi,\ \
M_\eta'(\varrho_2)=\psi^{-1}.
\]
Thus, for $\Delta_\eta\ne0$,
\[
\begin{array}{c|c}
|\psi|<1 & \varrho_1\text{ is attracting and }\varrho_2\text{ is repelling},\\
|\psi|>1 & \varrho_2\text{ is attracting and }\varrho_1\text{ is repelling},\\
|\psi|=1 & |w|\text{ is preserved under }w\mapsto\psi w,\quad M_\eta^m=\Id\iff\psi^m=1.
\end{array}
\]
When $\Delta_\eta=0$,
\[
\varrho_1=\varrho_2=\frac{\theta}{2}.
\]
For the attracting and repelling behavior of the roots and their multipliers $M_\eta'(\cdot)$, see \cite[Secs.~3.7.2--3.7.5]{Needham}.  If $c=0$ and $\theta\ne0$, then
\[
M_\eta(\delta)=\theta \ \ (\delta\ne0).
\]
If $c=\theta=0$, then the first iterate has $\delta=0$.

We next describe the fixed-point locus \(E_1\) of \(M_\eta\) and the equimodular set \(E_2\).  They will also be used to determine when
$\inf\limits_{n\ge0}\smin(D_n(z))>0$.
The cases \(uv=0\) and \(\Delta_\eta=0\) will be considered separately in the determination of the Julia set.

\noindent\textbf{The fixed-point hypersurface $E_1$}\ \
For $\delta\ne0$, $E_1$ is the fixed-point locus of $\widehat F$, and its projective closure in $\Pj^3$ is
\[
\begin{aligned}
E_1
&=
\left\{
[\theta:\delta:u:v]\in\Pj^3
\middle|
\delta^2-\theta\delta+r^2uv=0
\right\}\\
&=
\left\{
[z_0:z_1:z_2:z_3]\in\Pj^3
\middle|
(z_0-z_1)z_1=rz_2z_3
\right\},
\end{aligned}
\]
where $\theta-\delta=rz_1$.
Define
\[
E_{1,\mathrm{att}}
=\left\{
[\theta:\delta:u:v]\in E_1
\middle|
\delta\ne0,\ \left|\frac{r^2uv}{\delta^2}\right|<1
\right\},
\]
and
\[
E_{1,\mathrm{rep}}
=\left\{
[\theta:\delta:u:v]\in E_1
\middle|
\delta\ne0,\ \left|\frac{r^2uv}{\delta^2}\right|>1
\right\}.
\]

\begin{proposition}\label{p3}
Let $\xi\in\Pj^3$, and let
$\zeta=(\zeta_0,\zeta_1,\zeta_2,\zeta_3)\in\C^4\setminus\{0\}$ be a representative of $\xi$ with $\zeta_0\ne\zeta_1$.  Let
\[
\delta_\infty=\zeta_0-\zeta_1,\ \
a=\frac{\zeta_2}{\delta_\infty},\ \
b=\frac{\zeta_3}{\delta_\infty}.
\]
Then $\xi\in E_1$ if and only if
\begin{equation}\label{e41}
D_n(\zeta)=\delta_\infty(I+bX_n^*)(I+aX_n),\ \ 
\forall n\ge0.
\end{equation}
\end{proposition}

\begin{proof}
Using $X_n^*X_n=r(I+S_n)$, we obtain
\[
\delta_\infty(I+bX_n^*)(I+aX_n)
=\delta_\infty I+\zeta_2X_n+\zeta_3X_n^*
+\frac{r\zeta_2\zeta_3}{\delta_\infty}(I+S_n).
\]
Comparing with \eqref{e7}, the identity holds if and only if
\[
\delta_\infty\zeta_1=r\zeta_2\zeta_3.
\]
Since $\theta-\delta_\infty=r\zeta_1$, this is equivalent to
\[
\delta_\infty^2-\theta\delta_\infty+r^2\zeta_2\zeta_3=0,
\]
which is the defining equation of $E_1$ at $\xi$.
\end{proof}

For $\xi\in E_1$ with $\delta\ne0$, select a representative $\zeta$ as in Proposition~\ref{p3}.  Then
\[
M_\eta'(\delta_\infty)=\frac{r^2\zeta_2\zeta_3}{\delta_\infty^2}.
\]
Equivalently,
\[
E_{1,\mathrm{att}}
=\left\{\xi\in E_1\middle|\delta_\infty\ne0,\ |M_\eta'(\delta_\infty)|<1\right\},\ \
E_{1,\mathrm{rep}}
=\left\{\xi\in E_1\middle|\delta_\infty\ne0,\ |M_\eta'(\delta_\infty)|>1\right\}.
\]
By \eqref{e40}, hence
\[
E_{1,\mathrm{att}}\cap\Omega_{\mathrm{reg}}\subset\mathcal F(\widehat F).
\]  If
\[
[z]\in\mathcal F(\widehat F),\ \
\lim\limits_{k\to\infty}\widehat F^k([z])=\xi, \ \
\text{in }\Pj^3
\]
for some $\xi\in E_1\cap\Omega_{\mathrm{reg}}$, then $\xi\in E_{1,\mathrm{att}}$.

\begin{example}\label{x1}
Take
$\xi_r=[4r+1:1:2:2]$, $\zeta_r=(4r+1,1,2,2)$.
Then
$\delta=4r$, $\theta=5r$, $c=4r^2$,
so $\xi_r\in E_1$ and
$|M_\eta'(4r)|=\frac14<1$.
Thus $\xi_r\in E_{1,\mathrm{att}}$.  Since
$a=b=\frac1{2r}$,
from \eqref{e41}, then
\[
D_n(\zeta_r)=4r
\left(I+\frac1{2r}X_n^*\right)
\left(I+\frac1{2r}X_n\right).
\]
By $\|X_n\|\le r$, we have
\[
\left(I+\frac1{2r}X_n\right)^{-1}
=\sum\limits_{k=0}^{\infty}
\left(-\frac1{2r}X_n\right)^k,\ \
\left\|\left(I+\frac1{2r}X_n\right)^{-1}\right\|\le2.
\]
The same estimate holds for $I+(2r)^{-1}X_n^*$.  Hence
\[
\sup_{n\ge0}\|D_n(\zeta_r)^{-1}\|
\le\frac1r<\infty.
\]
Therefore Proposition~\ref{p1} implies
$\xi_r\notin p(D_{\pi,G})$.
 Moreover, it shows that $E_1$ is not entirely contained in the Julia set or the projective spectrum. For $r=2$, this result necessitates a revision of the corresponding inclusions in \cite[Theorem~B and Theorems~4.7--4.8]{ZYL}.
\end{example}

\noindent\textbf{The equimodular set $E_2$}\ \
For $[z]\in\Omega_{\mathrm{reg}}$, let $\varrho_1,\varrho_2$ be the two roots and let
$\psi=\frac{\varrho_2}{\varrho_1}$.
Then
$|\varrho_1|=|\varrho_2|\iff|\psi|=1$.
Since
$\frac{\theta^2}{4c}=\frac{(1+\psi)^2}{4\psi}$,
$c=r^2uv$,
if $\psi=e^{i\phi}$, then
$\frac{\theta^2}{4c}=\cos^2\frac\phi2$.
Let $x=\cos^2(\phi/2)$ and define
\begin{equation}\label{e42}
E_2
:=
\left\{
[\theta:\delta:u:v]\in\Pj^3
\middle|
\exists x\in[0,1],\ \theta^2=4r^2xuv
\right\}.
\end{equation}
On $\Omega_{\mathrm{reg}}$, \eqref{e42} can be written as
\begin{equation}\label{e43}
E_2\cap\Omega_{\mathrm{reg}}
=
\left\{[z]\in\Omega_{\mathrm{reg}}\middle|
|\varrho_1(z)|=|\varrho_2(z)|\right\}
=
\left\{[z]\in\Omega_{\mathrm{reg}}\middle|
|\psi([z])|=1\right\}.
\end{equation}
If $\psi=e^{2\pi i k/m}$ with $\gcd(k,m)=1$, then
$M_\eta^m=\Id$,
$\theta^2=4r^2uv\cos^2\frac{\pi k}{m}$.
Since $\mathbb{Q} \cap [0, 1]$ is dense in $[0, 1]$ and $x\mapsto\cos^2(\pi x)$ is continuous, the set
\[
\left\{\cos^2\frac{\pi k}{m}\middle|
 m\ge2,\ 1\le k<m,\ \gcd(k,m)=1\right\}
\]
is dense in $[0,1]$.  Hence
\[
E_2
=
\overline{
\left\{
 [z]\in\Omega_{\mathrm{reg}}\middle|
 \exists m\ge2,\ M_{\Pi([z])}^m=\Id
\right\}}.
\]
By \eqref{e17}, $\varrho_1$ and $\varrho_2$ are also the characteristic roots of the recurrence.

For $[z]\in E_2\cap\Omega_{\mathrm{reg}}$, let $\eta=\Pi([z])$.  Since $|\psi([z])|=1$, it follows from \eqref{e40} that the iterates of $\widehat F|_{L_\eta}$ form a normal family, so $[z]\in\mathcal F_\eta$.  The next proposition shows that, as a point of $\Pj^3$, $[z]$ belongs to $\mathcal J(\widehat F)$.
\begin{proposition}\label{p4}
\[
E_2\cap\Omega_{\mathrm{reg}}\subset\mathcal J(\widehat F).
\]
\end{proposition}

\begin{proof}
Let $[z_*] \in E_2 \cap \Omega_{\mathrm{reg}}$.  If $[z_*] \notin \Omega_\infty$, then $[z_*] \in \mathcal{J}(\widehat{F})$ by definition.  Assume $[z_*] \in \Omega_\infty$ and suppose $\{\widehat{F}^n\}$ is normal on an open neighborhood $U$ of $[z_*]$. In the chart $v=1$, we restrict to the slice $u=u_*$.  Since
\[
\frac{\partial\psi}{\partial\theta}
=-\frac{2\psi}{(\theta^2-4r^2u_*)^{1/2}}\ne0,
\]
$(\psi,w)$ are local holomorphic coordinates and
$\widehat F^n(\psi,w)=(\psi,\psi^nw)$.

Choose $q=(\psi_0,w_0)\in U$ satisfying $|\psi_0| = 1$ and $w_0 \neq 0$, and let
\[
p_n=\left(\psi_0\left(1+\frac1n\right),w_0\right).
\]
Then
$\lim\limits_{n\to\infty} p_n = q$, and
\[
\lim_{n\to\infty} \left| w(\widehat{F}^n(p_n)) - w(\widehat{F}^n(q)) \right|
= \lim_{n\to\infty} \left| \left( 1 + \frac{1}{n} \right)^n - 1 \right| |w_0|
= (e - 1)|w_0| > 0.
\]
If $\{\widehat{F}^m|_U\}_{m \ge 0}$ were locally equicontinuous at $q$, then $\lim\limits_{n\to\infty} p_n=q$ would imply
\[
\lim_{n\to\infty} \left| w(\widehat{F}^n(p_n)) - w(\widehat{F}^n(q)) \right| = 0,
\]
a contradiction.
	
Hence $\{\widehat{F}^n|_U\}_{n \ge 0}$ is not normal, and therefore $[z_*] \in \mathcal{J}(\widehat{F})$.
\end{proof}

We next relate $E_2$ to the zero sets of the polynomials $\mathsf G_n$ in $\Pj^3$.  Since $c=r^2uv$, the recursion \eqref{e17} shows by induction that $\mathsf G_n$ is a homogeneous polynomial of degree $n+1$ in $\theta,\delta,u,v$.  For $n\ge0$, define
\[
\mathcal H_n
:=
\left\{
[\theta:\delta:u:v]\in\Pj^3\middle|
\mathsf G_n(\theta,\delta,u,v)=0
\right\}.
\]
If $\delta_j(z)\ne0$ for $0\le j<n$, then \eqref{e33} implies
$[z]\in\mathcal H_n\iff\delta_n(z)=0$.
Since $\det A_{n+1}=\mathsf G_n$,
$[z]\in\mathcal H_n \iff\det A_{n+1}(z)=0$ for $z\in\Pj^3$.

Suppose now that $[z]\in\Omega_{\mathrm{reg}}$.  Since $\varrho_1\ne\varrho_2$, \eqref{e20} yields
\[
\mathsf G_n=
\frac{\varrho_1^{n+1}(\delta-\varrho_2)-
\varrho_2^{n+1}(\delta-\varrho_1)}
{\varrho_1-\varrho_2}.
\]
Consequently,
\begin{equation}\label{e44}
[z]\in\mathcal H_n
\iff
w([z])=\psi([z])^{-(n+1)}.
\end{equation}

\begin{figure}[H]
\centering
\paperfigure[width=0.96\textwidth]{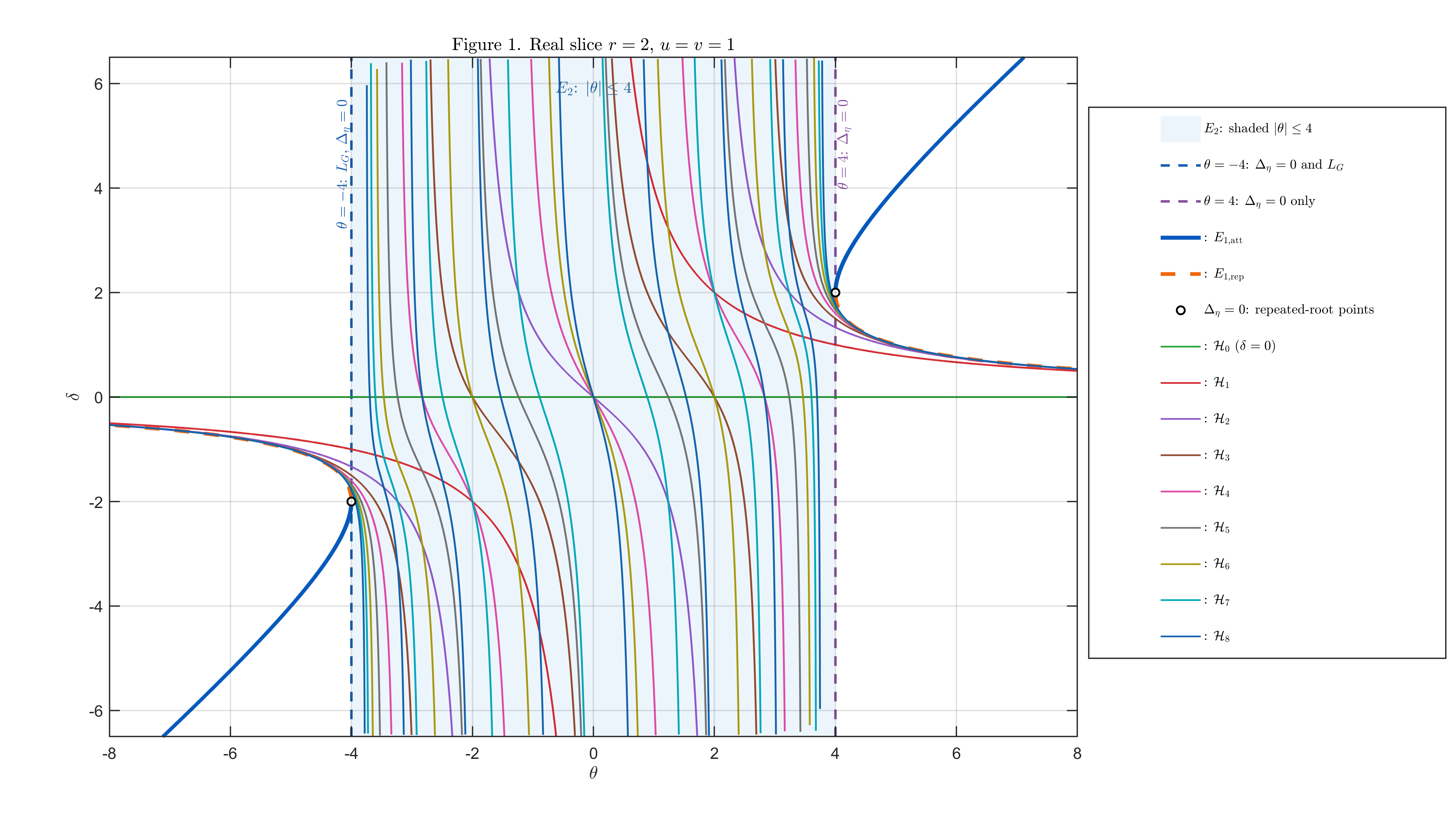}
\caption{The real affine slice $r=2$, $u=v=1$. The shaded strip $|\theta|\le4$ represents
$E_2$. The dashed lines $\theta=\pm4$ satisfy $\Delta_\eta=0$, and $\theta=-4$ also represents $L_G$. The dark-blue solid and orange dashed curves are $E_{1,\mathrm{att}}$ and $E_{1,\mathrm{rep}}$, respectively. The open circles mark the repeated-root points. The other coloured curves are $\mathcal H_n$, $0\le n\le8$, as indicated in the legend.}
\label{f1}
\end{figure}

Figure~\ref{f1} illustrates the relative position of these sets on a real affine slice.  The statements themselves are in $\Pj^3$.

Substituting $\delta_0=\mathsf G_0$ and
$\delta_j=\mathsf G_j/\mathsf G_{j-1}$ into \eqref{e30}, and then using polynomial continuation, implies
\begin{equation}\label{e45}
	\det D_{n+1}(z)
	=
	D_0(z)\,\mathsf G_n(z)^{r-1}
	\prod_{j=0}^{n-1}
	\mathsf G_j(z)^{(r-1)^2r^{\,n-j-1}}
	\qquad \forall n\ge0,
\end{equation}
with the empty product equal to $1$.

The next lemma gives the two accumulation results used in the proof of Theorem~\ref{t4}.

\begin{lemma}\label{l2}
On $\Omega_{\mathrm{reg}}$,
\begin{equation}\label{e46}
E_{1,\mathrm{rep}}
\subset
\overline{\bigcup_{n\ge0}\mathcal H_n},
\end{equation}
and
\begin{equation}\label{e47}
E_2\cap\Omega_{\mathrm{reg}}
\subset
\overline{\bigcup_{n\ge0}\mathcal H_n}.
\end{equation}
\end{lemma}

\begin{proof}
Without loss of generality, assume that $\varrho_2$ is the repelling fixed point, and let
$[z_*]=[\theta_*:\varrho_2:u_*:v_*]\in E_{1,\mathrm{rep}}$.
Then $|\psi_*|<1$.  Let
\[
w_n=\psi_*^{-(n+1)},\ \
\delta_n=\frac{w_n\varrho_2-\varrho_1}{w_n-1},\ \
[z_n]=[\theta_*:\delta_n:u_*:v_*].
\]
By \eqref{e44},
$[z_n]\in\mathcal H_n$.  Moreover, 
\[
\lim_{n\to\infty}|w_n|=\infty,\ \
\delta_n-\varrho_2=\frac{\varrho_2-\varrho_1}{w_n-1},\ \
\lim_{n\to\infty}\delta_n=\varrho_2.
\]
Hence $\lim\limits_{n\to\infty}[z_n]=[z_*]$ in $\Pj^3$.
Therefore \eqref{e46} holds.

Next let $[z_*]\in E_2\cap\Omega_{\mathrm{reg}}$.  Since $v_*\ne0$, take $v_*=1$ and let
$\psi_*=e^{i\phi}$, then $w_*=w(z_*)\in\Pj^1$.
For $w_*\in\C^*$ choose $h_*\in\C$ with $e^{h_*}=w_*$.  Define
\[
h_n=
\begin{cases}
h_*+\dfrac1{n+1},&w_*\in\C^*,\\[1mm]
-\sqrt{n+1},&w_*=0,\\[1mm]
\sqrt{n+1},&w_*=\infty,
\end{cases}\ \
w_n=e^{h_n}.
\]
Then
\[
\lim_{n\to\infty}w_n=w_*\ \ \text{in }\Pj^1,\ \
\lim_{n\to\infty}\frac{|h_n|}{n+1}=0.
\]
Choose $k_n\in\Z$ such that
\[
\left|\frac{2\pi k_n}{n+1}-\phi\right|\le\frac{\pi}{n+1},\ \
\psi_n=
\exp\!\left(\frac{2\pi i k_n-h_n}{n+1}\right).
\]
Hence
$\psi_n^{-(n+1)}=w_n$,
$\lim\limits_{n\to\infty}\psi_n=\psi_*$.
Moreover,
$\frac{\partial\psi}{\partial\theta}([z_*])\ne0$.
According to the holomorphic inverse function theorem \cite[Chap.~IV, Corollary~7.6]{Conway}, we obtain
\[
\exists\,(\theta_n)_{n\ge N},\ \
\lim_{n\to\infty}\theta_n=\theta_*,\ \
\forall n\ge N,\ \ \psi(\theta_n,u_*,1)=\psi_n.
\]
Since $\delta\mapsto w(\theta_n,\delta,u_*,1)$ is a M\"obius map, thus
\[
\exists\,(\delta_n)_{n\ge N},\ \
\lim_{n\to\infty}\delta_n=\delta_*,\ \
\forall n\ge N,\ \ w(\theta_n,\delta_n,u_*,1)=w_n.
\]
Let $z_n=(\theta_n,\delta_n,u_*,1)$. By \eqref{e44}, we have
\[
\lim_{n\to\infty}[z_n]=[z_*],\ \
w([z_n])=\psi([z_n])^{-(n+1)},\ \
[z_n]\in\mathcal H_n,\quad n\ge N.
\]
Therefore \eqref{e47} holds.
\end{proof}

\begin{figure}[H]
\centering
\paperfigure[width=0.99\textwidth]{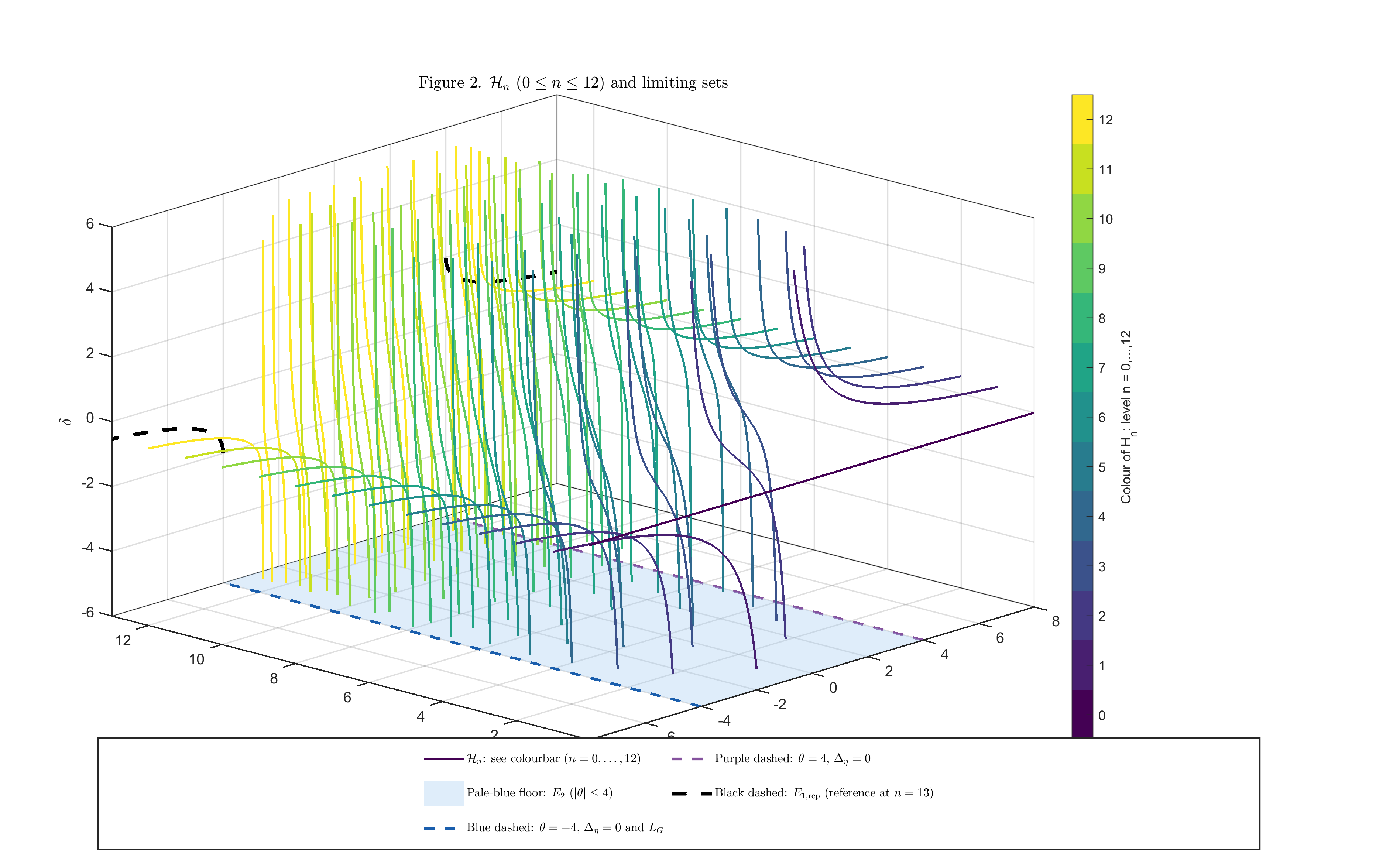}
\caption{The curves $\mathcal H_n$, $0\le n\le12$, on the slice $r=2$, $u=v=1$, stacked at auxiliary level $n$, with vertical coordinate $\delta$. Colours indicate $n$ as shown by the colourbar. The shaded floor represents $E_2$, bounded by $\theta=\pm4$. the left boundary also represents $L_G$. The black dashed curve at level $n=13$ indicates $E_{1,\mathrm{rep}}$. The projections onto the $(\theta,\delta)$-plane illustrate accumulation toward $E_{1,\mathrm{rep}}$ and $E_2$.}
\label{f2}
\end{figure}

Figure~\ref{f2} illustrates the accumulation described in Lemma~3.8 by projecting the stacked curves $\mathcal H_n$ onto the $(\theta,\delta)$-plane.

\begin{theorem}\label{t4}
The Julia set is the closure of the finite singularity hypersurfaces, namely
\[
\mathcal J(\widehat F)=\overline{\bigcup_{n\ge0}\mathcal H_n}\subset p(D_{\pi,G}).
\]
On $\Omega_{\mathrm{reg}}$,
\[
\mathcal J(\widehat F)\cap\Omega_{\mathrm{reg}}
=
\left[
\left(\bigcup_{n\ge0}\mathcal H_n\right)
\cup E_{1,\mathrm{rep}}
\cup E_2
\right]\cap\Omega_{\mathrm{reg}},
\]
whereas
\[
E_{1,\mathrm{att}}\subset\mathcal F(\widehat F).
\]
\end{theorem}

\begin{proof}
Let
$\mathcal C=\overline{\bigcup\limits_{n\ge0}\mathcal H_n}$.
We first show that $\mathcal J(\widehat F)\subset\mathcal C$.  Let $[z_*]\notin\mathcal C$.
	
If $[z_*]\in\Omega_{\mathrm{reg}}$, by
Lemma~\ref{l2}, then
$[z_*]\notin E_{1,\mathrm{rep}}\cup E_2$. Hence
$|\varrho_1(z_*)|\ne|\varrho_2(z_*)|$. Without loss of generality,
assume $|\varrho_2(z_*)|<|\varrho_1(z_*)|$.  Then $|\psi([z_*])|<1$. Since $[z_*]\notin E_{1,\mathrm{rep}}$, we have
$\delta_*\ne\varrho_2(z_*)$, and hence $w([z_*])\in\C$. Choose a
neighborhood $U$ of $[z_*]$ with $U\cap\mathcal C=\varnothing$ and
$|\psi([z])|\le q<1$, $|w([z])|\le M$, $[z]\in U$.
Thus
\[
\lim_{n\to\infty}\sup_{[z]\in U}|\psi([z])^n w([z])|=0.
\]
Hence $U\subset\Omega_\infty$ and $\{\widehat F^n|_U\}_{n\ge0}$ is normal.
	
If $[z_*]\in(\Pj^3\setminus\Omega_{\mathrm{reg}})\setminus\mathcal C$,
then $u_*v_*\ne0$ or $u_*v_*=0$ with $\theta_*\ne0$. In the first case,
$\Delta_\eta([z_*])=0$ and
$[z_*]\in\overline{E_2\cap\Omega_{\mathrm{reg}}}\subset\mathcal C$, a contradiction. In the second case, use the affine chart $\theta=1$. At
$[z_*]$ the two roots are $1$ and $0$. Since $[z_*]\notin\mathcal H_0$, we
have $\delta_*\ne0$. After shrinking $U$ so that
$U\cap\mathcal C=\varnothing$, choose the roots holomorphically on $U$
and arrange that
\[
\left|\frac{\varrho_2(z)}{\varrho_1(z)}\right|\le q<1,\ \
\left|\frac{\delta(z)-\varrho_1(z)}{\delta(z)-\varrho_2(z)}\right|\le M,
\ \ z\in U.
\]
Then
\[
\lim_{n\to\infty}\sup_{[z]\in U}|\psi([z])^n w([z])|=0,
\]
so $U\subset\Omega_\infty$ and $\{\widehat F^n|_U\}_{n\ge0}$ is
normal by \eqref{e40}. The remaining case
$u_*v_*=0$, $\theta_*=0$ is impossible because $\mathsf G_1(z_*)=0$.
Therefore $[z_*]\notin\mathcal J(\widehat F)$, so
$\mathcal J(\widehat F)\subset\mathcal C$.
	
Conversely, let $\xi=[z]\in\mathcal H_n$ and let $m\le n$ be minimal
with $\mathsf G_m(z)=0$. Then
$\delta_j(z)\ne0, 0\le j<m$, $\delta_m(z)=0$.
If $uv=0$, the $m$-th iterate lies in $I(\widehat F)$, otherwise its next image is
$[0:1:0:0]\in I(\widehat F)$. Thus $\xi\notin\Omega_\infty$, and hence $\xi\in\mathcal J(\widehat F)$. Therefore,
$\mathcal H_n\subset\mathcal J(\widehat F)$ for every $n\ge0$. Since
$\mathcal J(\widehat F)$ is closed,
$\mathcal C\subset\mathcal J(\widehat F)$.
Thus $\mathcal J(\widehat F)=\mathcal C$. By
\eqref{e45} and
Proposition~\ref{p1},
$\mathcal H_n\subset p(D_{\pi,G})$ for every $n\ge0$. Since $p(D_{\pi,G})$ is
closed, hence
\[
\mathcal J(\widehat F)\subset p(D_{\pi,G}).
\]
	
On $\Omega_{\mathrm{reg}}$, the stated description of
$\mathcal J(\widehat F)$ follows from \eqref{e44},
Lemma~\ref{l2}, and
Proposition~\ref{p4}. By
\eqref{e40},
$E_{1,\mathrm{att}}\cap\Omega_{\mathrm{reg}}\subset\mathcal F(\widehat F)$.
For $\xi\in E_{1,\mathrm{att}}\setminus\Omega_{\mathrm{reg}}$,
\eqref{e38} and $|r^2uv/\delta^2|<1$ imply $uv=0$ and
$\delta=\theta\ne0$. Choose the roots on a neighborhood $U$ of $\xi$
so that
\[
\left|\frac{\varrho_2}{\varrho_1}\right|\le q<1,\ \
|w|\le M.
\]
Then \eqref{e40} implies
$\lim\limits_{n\to\infty}\sup\limits_{[z]\in U}|\psi([z])^n w([z])|=0$.
Hence $E_{1,\mathrm{att}}\subset\mathcal F(\widehat F)$.
\end{proof}

\begin{figure}[H]
\centering
\paperfigure[width=0.79\textwidth]{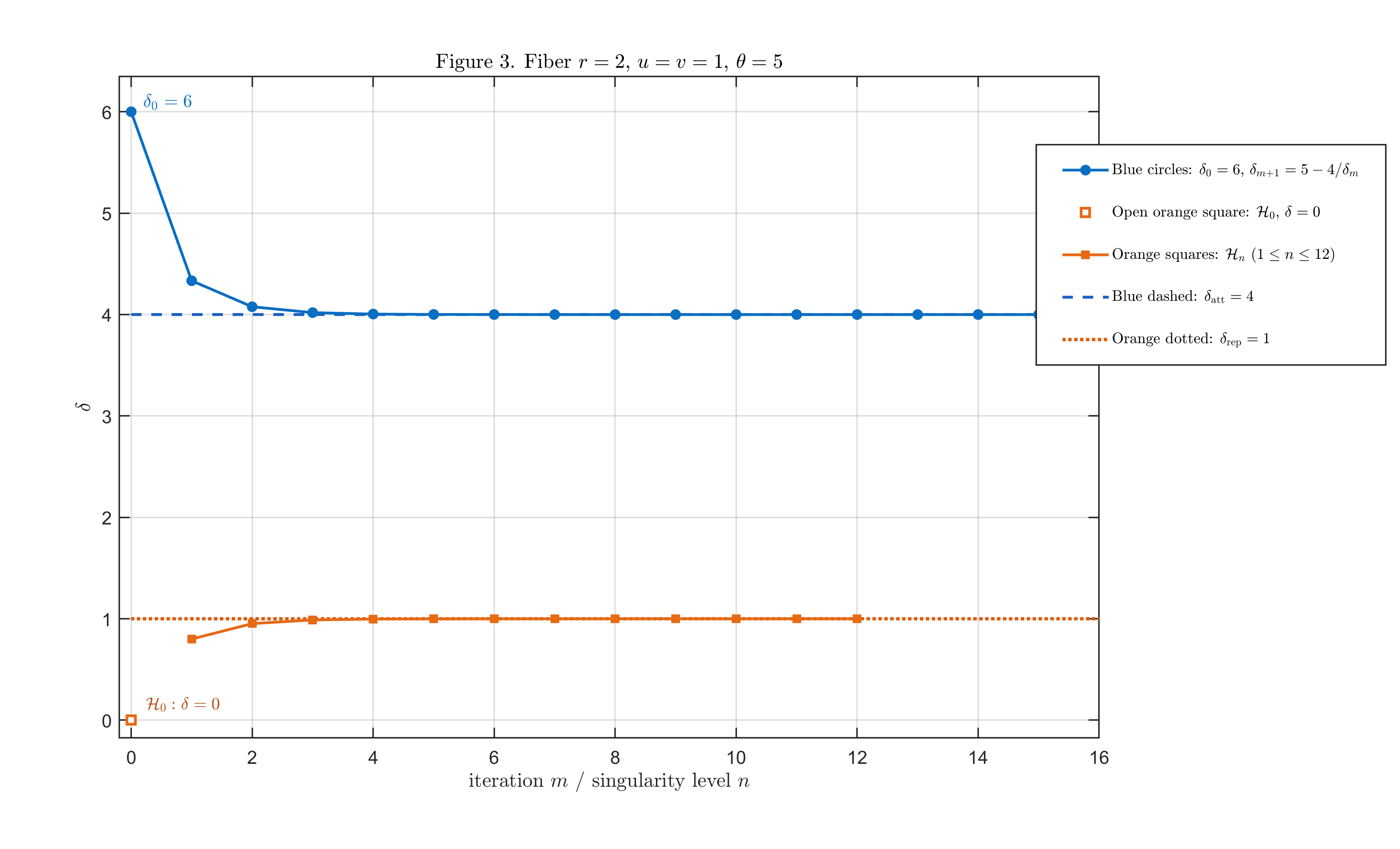}
\caption{The fiber $r=2$, $u=v=1$, $\theta=5$. The blue orbit starting at $\delta_0=6$ converges to the attracting fixed point $\delta_{\mathrm{att}}=4$. The orange squares represent $\mathcal H_n$, $0\le n\le12$, and converge to the repelling fixed point $\delta_{\mathrm{rep}}=1$. The open square marks $\mathcal H_0$ at $\delta=0$.}
\label{f3}
\end{figure}

Figure~\ref{f3} shows the two limiting behaviors in Theorem~\ref{t4} on a single real fiber.
\FloatBarrier

\noindent\textbf{Unbounded inverse norms}\ \
Let $[\zeta]\in E_{1,\mathrm{att}}$ be a fixed point of $\widehat F$, where $\zeta\in\mathbb C^4\setminus\{0\}$ is a representative.  By \eqref{e41}, we have
\[
D_n(\zeta)=\delta_\infty(I+bX_n^*)(I+aX_n)=\delta_\infty(I+brY_n^*)(I+arY_n).
\]
It suffices to study $I+aX_n$.  We determine for which $a\in\C$
\[
\forall n\ge0,\ \ I+aX_n\text{ invertible},\ \
\sup_{n\ge0}\bigl\|(I+aX_n)^{-1}\bigr\|<\infty.
\]
By Lemma~\ref{l1}, $Y_n^k$ is a partial isometry, for $\forall k\ge0$.
Next, we first compute the rank of $Y_n^k$ and then decompose $Y_n|_{\mathscr H_n^0}$.

\begin{lemma}\label{l3}
For $n\ge0$, one has
\[
\operatorname{rank}(Y_n^k)=r^{n-k}\ \ (0\le k\le n),\ \
\forall k\ge n,\ \ \operatorname{rank}(Y_n^k)=1.
\]
\end{lemma}

\begin{proof}
For $k\ge1$,
$G_0G_1\cdots G_{k-1}
=\prod_{j=0}^{k-1}G_j
\cong G^k$.
By Lemma~\ref{l1}, we have
$$E_k=\frac1{|G|^k}\sum\limits_{h\in G^k}h,$$
hence $\pi_n(E_k)$ is the orthogonal projection onto
$\mathscr H_n^{G^k}=\{f\in\ell^2(G^n)\mid f(hx)=f(x)\ \forall h\in G^k\}$.
Moreover,
\[
\mathscr H_n^{G^k}
=\bigoplus_{\mathcal O\in G^k\backslash G^n}\C\mathbf 1_{\mathcal O},\ \
\operatorname{rank}\pi_n(E_k)
=\dim \mathscr H_n^{G^k}
=|G^k\backslash G^n|,
\]
where $\mathbf 1_{\mathcal O}$ is the characteristic function of the orbit $\mathcal O$.

For $k \le n$, every $h \in G^k$ acts only on the first $k$ coordinates. If $h$ fixes a word $x \in G^n$, then it fixes every point of the cylinder $xG^{\mathbb N}$. Since $G$ is finite abelian, no nonidentity element can fix a nonempty cylinder, see \cite[p.~203]{KSS}. Hence the action of $G^k$ on $G^n$ is free. Thus
\[
\operatorname{Stab}_{G^k}(x)=\{1\},\ \ x\in G^n,\ \
|G^k\backslash G^n|
=\frac{|G^n|}{|G^k|}=r^{n-k}.
\]
If $k\ge n$, the action of $G^n$ on $G^n$ is free and $|G^n|=r^n$, hence transitive.  Therefore $G^k$ also has one orbit.

Finally,
$Y_n^k=\pi_n(t^{-k})\pi_n(E_k)$,
and $\pi_n(t^{-k})$ is unitary.  Therefore $\operatorname{rank}(Y_n^k)=\operatorname{rank}(\pi_n(E_k))$.  This proves the lemma.
\end{proof}

Let
$\mathscr H_n^0=\mathbf 1_n^\perp$,
$N_n=Y_n|_{\mathscr H_n^0}$.
Then $N_n^n=0$ and
$\operatorname{rank}(N_n^k)=r^{n-k}-1$,
$0\le k\le n$.

\begin{theorem}\label{t5}
For every $n\ge1$, let $J_m$ be the nilpotent backward shift. Then
\begin{equation}\label{e48}
Y_n\cong
1\oplus J_n^{\oplus(r-1)}
\oplus\bigoplus_{m=1}^{n-1}
J_m^{\oplus (r-1)^2r^{n-m-1}}.
\end{equation}
\end{theorem}

\begin{proof}
By Lemma~\ref{l3}, $N_n^n=0$ and $N_n^k$ is a partial isometry for $k\ge0$.  By the Halmos--Wallen decomposition theorem \cite[Theorem~1]{HalmosWallen}, there exist integers $d_1,\ldots,d_n\ge0$ such that
\[
N_n\cong\bigoplus_{m=1}^{n}J_m^{\oplus d_m}.
\]
For $1\le m\le n$,
\[
\sum_{j=m}^{n}d_j
=\operatorname{rank}(N_n^{m-1})-\operatorname{rank}(N_n^m).
\]
By Lemma~\ref{l3},
$\sum\limits_{j=m}^{n}d_j=(r-1)r^{n-m}$.
It follows that
\[
\begin{aligned}
d_m
&=(r-1)r^{n-m}-(r-1)r^{n-m-1}
=(r-1)^2r^{n-m-1},
&&1\le m<n,\\
d_n&=r-1.
\end{aligned}
\]
For $1\le m\le n$, select orthonormal vectors
$e_{m,l,0},\ldots,e_{m,l,m-1}$, $1\le l\le d_m$,
such that
\[
N_ne_{m,l,0}=0,\ \
N_ne_{m,l,j}=e_{m,l,j-1}
\quad(1\le j<m).
\]
Thus
\[
\mathscr H_n
=
\C\mathbf 1_n
\oplus
\bigoplus_{m=1}^{n}\bigoplus_{l=1}^{d_m}
\operatorname{span}\{e_{m,l,0},\ldots,e_{m,l,m-1}\}.
\]
Let $f_{m,l,0},\ldots,f_{m,l,m-1}$ be the canonical basis of the $(m,l)$--th copy of $\C^m$ and define
\[
U_n:\mathscr H_n\longrightarrow
\C\oplus\bigoplus_{m=1}^{n}(\C^m)^{\oplus d_m},\ \
U_n\frac{\mathbf 1_n}{\|\mathbf 1_n\|}=1,\ \
U_ne_{m,l,j}=f_{m,l,j}.
\]
The map $U_n$ is unitary and
\[
U_nY_nU_n^*
=1\oplus\bigoplus_{m=1}^{n}J_m^{\oplus d_m}.
\]
Hence \eqref{e48} follows.
\end{proof}

By \eqref{e48} and \eqref{e15}, we obtain
\begin{equation}\label{e49}
D_n(z)\cong
D_0(z)
\oplus
\bigl(\mathcal{A}_n(z)\bigr)^{\oplus(r-1)}
\oplus
\bigoplus_{m=1}^{n-1}
\bigl(\mathcal{A}_m(z)\bigr)^{\oplus (r-1)^2r^{n-m-1}}.
\end{equation}

To continue, we define
\begin{equation}\label{e50}
\mathcal F_\infty
=
\left\{
[z]\in\mathcal F(\widehat F)
\middle|
\begin{array}{c}
\exists\,\xi=[\xi_0:\xi_1:\xi_2:\xi_3]\in E_{1,\mathrm{att}}\text{ such that }
\displaystyle\lim\limits_{k\to\infty}\widehat F^k([z])=\xi\text{ in }\Pj^3,\\[1mm]
\displaystyle
\max\left\{
\left|\frac{\xi_2}{\xi_0-\xi_1}\right|,
\left|\frac{\xi_3}{\xi_0-\xi_1}\right|
\right\}\ge\frac1r
\end{array}
\right\}.
\end{equation}
For the fixed point $\xi$, the ratios $a$ and $b$ are defined as in Proposition~\ref{p3}, they do not depend on the chosen representative.  

For a finite-dimensional matrix $T$, its smallest singular value is
\begin{equation}\label{e51}
\smin(T)
=\min_{\norm{x}=1}\norm{Tx}
=\sqrt{\lambda_{\min}(T^*T)}.
\end{equation}
If $T$ is invertible, then
\begin{equation}\label{e52}
\smin(T)=\norm{T^{-1}}^{-1},
\end{equation}
and for matrices of the same size,
\begin{equation}\label{e53}
\bigl|\smin(T)-\smin(S)\bigr|\le\norm{T-S}.
\end{equation}

\begin{lemma}\label{l4}
Fix $r$ and $z$ with $\delta\ne0$.  If $D_n(F(z))$ is invertible, then $D_{n+1}(z)$ is invertible and
\begin{equation}\label{e54}
\norm{D_n(F(z))^{-1}}
\le
\norm{D_{n+1}(z)^{-1}}
\le
C(z)\bigl(1+\norm{D_n(F(z))^{-1}}\bigr),
\end{equation}
where $C(z)>0$ is independent of $n$.
\end{lemma}

\begin{proof}
By Proposition~\ref{p2}, we have
\[
D_{n+1}(z)^{-1}
=
\begin{pmatrix}I&0\\-(\delta Q_n)^{-1}ruQ_n\mathcal U_n^*P_n&I\end{pmatrix}
\begin{pmatrix}\mathcal T_n^{-1}&0\\0&(\delta Q_n)^{-1}\end{pmatrix}
\begin{pmatrix}I&-rvP_n\mathcal U_nQ_n(\delta Q_n)^{-1}\\0&I\end{pmatrix}
\]
and $\mathcal T_n$ is unitarily equivalent to $D_n(F(z))$, hence $\norm{D_n(F(z))^{-1}}\le\norm{D_{n+1}(z)^{-1}}$.

Moreover,
$\|P_n\|=\|Q_n\|=\|\mathcal U_n\|=1$,
$\norm{rvP_n\mathcal U_nQ_n(\delta Q_n)^{-1}}\le\frac{r|v|}{|\delta|}$,
$\norm{(\delta Q_n)^{-1}ruQ_n\mathcal U_n^*P_n}\le\frac{r|u|}{|\delta|}$.
Therefore
\[
\norm{D_{n+1}(z)^{-1}}
\le
\left(1+\frac{r|u|}{|\delta|}\right)
\left(1+\frac{r|v|}{|\delta|}\right)
\max\left\{\norm{D_n(F(z))^{-1}},|\delta|^{-1}\right\}.
\]
Taking
$C(z)=\left(1+\frac{r|u|}{|\delta|}\right)\left(1+\frac{r|v|}{|\delta|}\right)\max\{1,|\delta|^{-1}\}$, then we obtain \eqref{e54}.
\end{proof}

\begin{proposition}\label{p5}
We prove that,
\[
\mathcal F_\infty\subset \mathcal F(\widehat F)\cap p(D_{\pi,G}).
\]
More precisely, if $[z]\in \mathcal F_\infty$, then every homogeneous representative $z$ satisfies
\[
\inf_{n\ge0}\smin(D_n(z))=0.
\]
\end{proposition}

\begin{proof}
We first prove
\begin{equation}\label{e55}
\sup_{n\ge0}\norm{(I+aX_n)^{-1}}<\infty
\iff |a|<\frac1r.
\end{equation}
If $|a|<1/r$, the Neumann series implies
\[
\forall n\ge0,\ \
\norm{(I+aX_n)^{-1}}
\le \sum_{j=0}^{\infty}(r|a|)^j
=\frac{1}{1-r|a|}.
\]
If $|a|\ge1/r$, by Theorem~\ref{t5}, whenever $I+aX_n$ is invertible, one has
\[
\norm{(I+aX_n)^{-1}}
\ge
\norm{(I+arJ_n)^{-1}e_{n-1}}
=
\left(\sum_{k=0}^{n-1}|ar|^{2k}\right)^{1/2}.
\]
Moreover,
$\lim\limits_{n\to\infty}\sum\limits_{k=0}^{n-1}|ar|^{2k}=\infty$.
This proves \eqref{e55}.  Taking adjoints gives the same statement for $I+bX_n^*$. 

Let $\xi$ be the attracting fixed point in \eqref{e50}. By Remark~\ref{a1}, let $\zeta \in \C^4\setminus \{0\}$ represent $\xi$. Assume first that $|a|\ge1/r$, $|b|<1/r$.  From Theorem~\ref{t5} and \eqref{e52}, we have
\[
0\le \smin(I+aX_n)\le \smin(I+arJ_n)
=\norm{(I+arJ_n)^{-1}}^{-1},\ \
\lim\limits_{n\to\infty}\smin(I+aX_n)=0.
\]
Thus \eqref{e41} implies $\lim\limits_{n\to\infty}\smin(D_n(\zeta))=0$.  The same conclusion holds when $|b|\ge1/r$ and $|a|<1/r$. 

Select a nonzero representative $z\in\mathbb C^4$ of $[z]$ and let
$z_k=F^k(z)$.  Since $F$ leaves $\theta,u,v$ unchanged and $\lim\limits_{k\to\infty}\widehat F^k([z])=\xi$ in $\Pj^3$, then Remark~\ref{a1} yields $\lim\limits_{k\to\infty}z_k=\zeta$ in $\C^4$.

Choose $k_n$ so that $\norm{D_n(z_{k_n})-D_n(\zeta)}<\frac1n$, $k_n\ge n$, $\forall n\ge1$.
By \eqref{e53},
\[
0\le \smin(D_n(z_{k_n}))\le \smin(D_n(\zeta))+\frac1n,\ \
\lim\limits_{n\to\infty}\smin(D_n(z_{k_n}))=0.
\]
If $D_n(z_{k_n})$ is singular for some $n$, Proposition~\ref{p2} shows that $D_{n+k_n}(z)$ is singular, and hence
\[
\inf_{m\ge0}\smin(D_m(z))=0.
\]
If all $D_n(z_{k_n})$ are invertible, then it follows from \eqref{e52} and Lemma~\ref{l4} that
\[
\norm{D_{n+k_n}(z)^{-1}}\ge \smin(D_n(z_{k_n}))^{-1},\ \
\lim\limits_{n\to\infty}\norm{D_{n+k_n}(z)^{-1}}=\infty,
\]
thus $\inf\limits_{m\ge0}\smin(D_m(z))=0$. Therefore, from Proposition~\ref{p1} we obtain $\mathcal F_\infty\subset \mathcal F(\widehat F)\cap p(D_{\pi,G})$.
\end{proof}

\begin{proposition}\label{p5.1}
Let $L_G$ be as in \eqref{e34}, then $L_G \subset \mathcal J(\widehat F)\cup \mathcal F_\infty$.
\end{proposition}
\begin{proof}
For $\forall [z] \in L_G\setminus\mathcal J(\widehat F)$, we have 
$\lim\limits_{n \to \infty}\widehat F^n([z])=\xi\in E_{1,\mathrm{att}}$.
From Proposition~\ref{p3} and taking $n=0$, hence 
$D_0(\zeta)=\delta_\infty(1+ra)(1+rb)=0$. Thus $a=-\frac1r$ or $b=-\frac1r$ and then 
$[z]\in\mathcal F_\infty$ by Proposition~\ref{p5}. Therefore
$L_G \subset \mathcal J(\widehat F)\cup \mathcal F_\infty$.
\end{proof}
In Figure~\ref{f1}, $L_G$ is precisely the left vertical dashed line $\theta=-4$.

\begin{theorem}\label{t6}
For every nontrivial finite abelian group $G$, one has
\[
p(D_{\pi,G})=\mathcal J(\widehat F)\cup \mathcal F_\infty.
\]
\end{theorem}

\begin{proof}
By \eqref{e37}, Theorem~\ref{t4}, Proposition~\ref{p5}, and Proposition~\ref{p5.1}, we have
$\mathcal J(\widehat F)\cup \mathcal F_\infty\subset p(D_{\pi,G})$.

Conversely, let $[z]\notin \mathcal J(\widehat F)\cup \mathcal F_\infty$.  By \eqref{e45}, $D_n(z)$ is invertible for every $n\ge0$.
There is $\xi\in E_{1,\mathrm{att}}$ with $|a|,|b|<1/r$ such that $\lim\limits_{k\to\infty}\widehat F^k([z])=\xi$ in $\Pj^3$.  Let $z_k=F^k(z)$, and choose the normalized representative $\zeta$ of $\xi$ as in Remark~\ref{a1}. Then $\lim\limits_{k\to\infty}z_k=\zeta$ in $\C^4$.

For $w, w' \in \C^4$, since $\norm{S_n}\le r-1$ and $\norm{X_n}\le r$, \eqref{e7} implies
\[
\begin{aligned}
\norm{D_n(w)-D_n(w')}
&\le |w_0-w_0'|+(r-1)|w_1-w_1'|
+r|w_2-w_2'|+r|w_3-w_3'|\\
&\le \sqrt{1+(r-1)^2+2r^2}\,\norm{w-w'}_2.
\end{aligned}
\]
Let
\[
C_r:=\sqrt{1+(r-1)^2+2r^2},\ \
M:=\sup\limits_{n\ge0}\norm{D_n(\zeta)^{-1}}<\infty
\]
Then
\[
\sup_{n\ge0}\norm{D_n(w)-D_n(w')}
\le C_r\norm{w-w'}_2.
\]
Select $k_*$ so that $MC_r\norm{z_{k_*}-\zeta}_2<1/2$, and let $E_n=D_n(\zeta)^{-1}(D_n(z_{k_*})-D_n(\zeta))$.  Then
$\norm{E_n}<1/2$ for $\forall n\ge0$, and hence
\[
D_n(z_{k_*})^{-1}=(I+E_n)^{-1}D_n(\zeta)^{-1},\ \
\sup_{n\ge0}\norm{D_n(z_{k_*})^{-1}}\le2M.
\]
Let $M_j:=\sup\limits_{n\ge0}\norm{D_n(z_j)^{-1}}$.  Lemma~\ref{l4} implies $M_j\le\max\{|D_0(z)|^{-1},C(z_j)(1+M_{j+1})\}<\infty$ for $j=k_*-1,\ldots,0$.  Thus $M_0<\infty$, and Proposition~\ref{p1} yields $[z]\notin p(D_{\pi,G})$.
\end{proof}

\begin{example}\label{x2}
For $G=\Z_2$, the projective spectrum proposed in \cite[Conjecture~4.9]{ZYL} requires the additional term $\mathcal F_\infty$.  For every $r\ge2$, consider
\begin{equation}\label{e56}
z^{(r)}=\left(3-\frac1r:-\frac1r:\frac4r:\frac1{5r}\right) \in \C^4.
\end{equation}
For this point, $\theta=2$, $\delta_0=3$, and $c=r^2uv=4/5$.  The roots of
\[
\delta^2-2\delta+\frac45=0
\]
are $\delta_\pm=1\pm1/\sqrt5$.  Since
\[
\left|\frac{4/5}{\delta_+^2}\right|<1,\ \
\left|\frac{u}{\delta_+}\right|=\frac4{r\delta_+}>\frac1r,
\]
we have $[z^{(r)}]\in \mathcal F_\infty\subset\mathcal F(\widehat F)\cap p(D_{\pi,G})$.

Also, $D_0(z^{(r)})=31/5\ne0$ and $\theta^2/(4c)=5/4>1$.  Let $\gamma=2/\sqrt5$, $x=\sqrt5/2$, and select $s>0$ with $x=\cosh s$.  The Chebyshev identity in Subsection~2.3 implies
\[
0<\gamma\frac{U_{n-1}(x)}{U_n(x)}
=\gamma\frac{\sinh(ns)}{\sinh((n+1)s)}<\gamma<1<\delta_0.
\]
Hence $\mathsf G_n(z^{(r)})\ne0$ for $n\ge0$, and the initial point is not the repelling fixed point.  Thus
\[
[z^{(r)}]\in \mathcal F_\infty\setminus\bigl(L_G\cup\mathcal J(\widehat F)\bigr).
\]

On the slice $z^{(r)}_1=-1/r$, $z^{(r)}_2=4/r$, $z^{(r)}_3=1/(5r)$, let $z^{(r)}_0$ vary.  Then $\theta=z_0-(r-1)/r$ and $c=4/5$.  For $[z^{(r)}]\notin\mathcal J(\widehat F)$ with $|\varrho_1|\ne|\varrho_2|$, let $\varrho$ be the root of larger modulus.  Since $\theta=\varrho+c/\varrho$, Theorem~\ref{t6} and \eqref{e50} show that $D_{\pi,G}(z^{(r)})$ is invertible exactly when $|\varrho|>4$.  Setting $\varrho=4e^{it}$ on the boundary yields
\[
z^{(r)}_0(t)=\frac{r-1}{r}+4e^{it}+\frac15e^{-it}
=\frac{r-1}{r}+\frac{21}{5}\cos t+i\frac{19}{5}\sin t.
\]
By \eqref{e15}, the matrix $\mathcal A_m-\theta I_m$ is diagonally similar to a real symmetric tridiagonal matrix. Hence the singularities of $\mathcal A_m$ are real and bounded by
\[
|\theta| \le 1+\frac{4}{\sqrt{5}} < \frac{21}{5}.
\]
The remaining Julia contributions in Theorem 3.9 are contained in
\[
\left[-\frac{4}{\sqrt{5}}, \frac{4}{\sqrt{5}}\right] \cup \left\{-\frac{9}{5}\right\}.
\]
Consequently, all Julia contributions lie inside the ellipse.
Therefore, the set of $z^{(r)}_0$ for which the pencil is not invertible forms the closed ellipse
\[
\frac{\bigl(\Re z^{(r)}_0-(r-1)/r\bigr)^2}{(21/5)^2}
+\frac{(\Im z^{(r)}_0)^2}{(19/5)^2}\le1,
\]
with semiaxes $21/5$ and $19/5$.  Equivalently, the spectrum of
\[
A_r=\pi\!\left(\frac1rC_G-4Y-\frac15Y^*\right)
\in C^*(\pi(\Gamma_G))
\]
is this closed ellipse.  In particular, $\sigma(A_r)$ has positive planar Lebesgue measure, with area $399\pi/25$.

Figure~\ref{f4} illustrates the spectral ellipse for $r=2$. It also shows its boundary lifted to the auxiliary height $t/(2\pi)$, along with the sample point $z^{(2)}$ from \eqref{e56}.

\begin{figure}[H]
\centering
\paperfigure[width=0.85\textwidth]{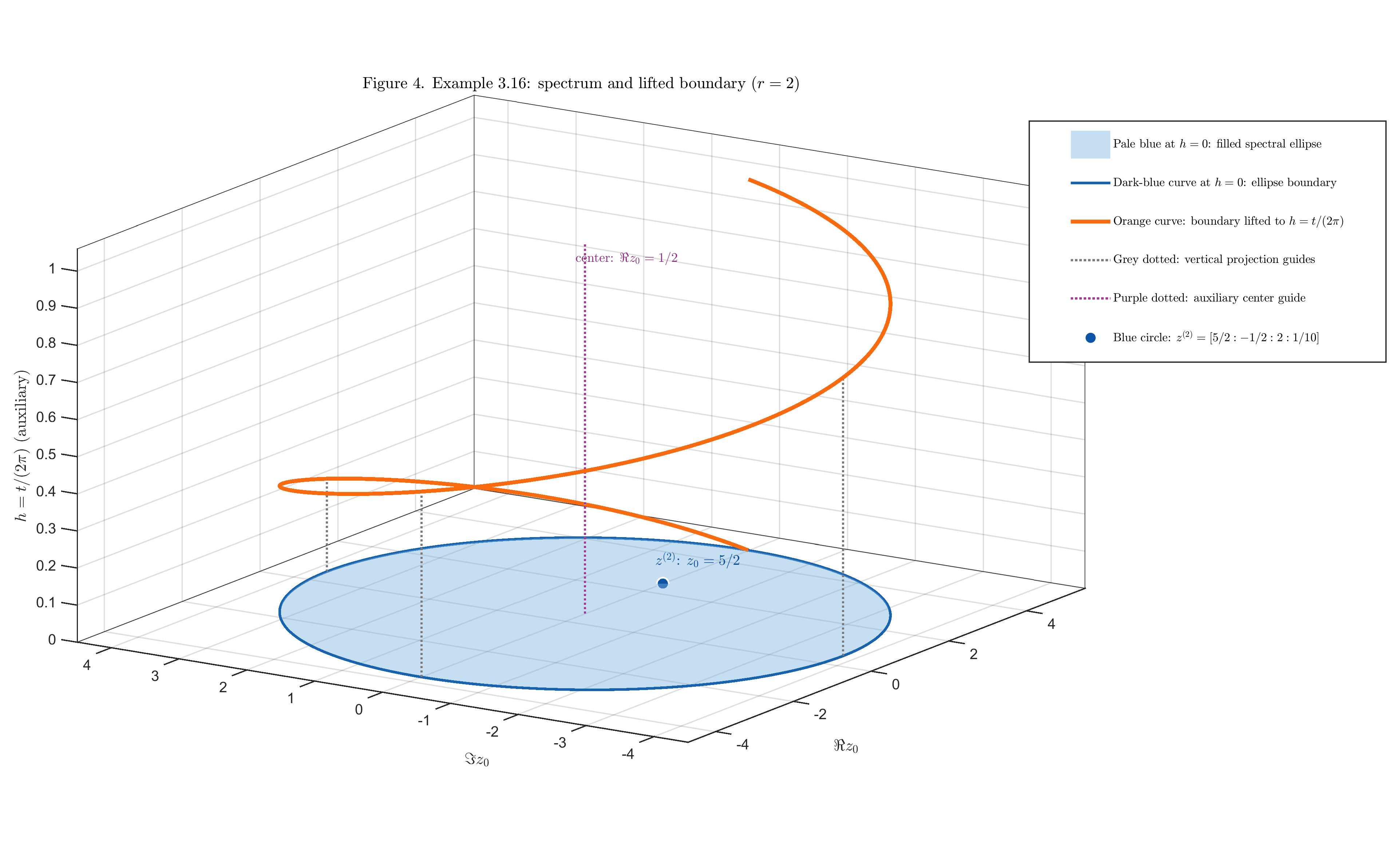}
\caption{Example~3.16 for $r=2$. The shaded ellipse at height zero represents the spectrum, and the dark-blue curve is its boundary. The orange curve is the boundary lifted
to the auxiliary height $t/(2\pi)$, not an additional part of the spectrum. The marked point is $[5/2:-1/2:2:1/10]\in \mathcal F_\infty\setminus (L_G\cup\mathcal J(\widehat F))$.}
\label{f4}
\end{figure}
\end{example}


\subsection{Weak equivalence and equality of the projective spectra}\label{ss2}

For unitary representations $\pi_1$ and $\pi_2$ of $\Gamma_G$, the notation $\pi_1\prec\pi_2$ stands for weak containment, and $\pi_1\sim_w\pi_2$ for weak equivalence. Equivalently,
\[
\pi_1\sim_w\pi_2
\iff
\norm{\pi_1(a)}=\norm{\pi_2(a)},\ \
\text{for every }a\in\C\Gamma_G.
\]
The reader is referred to \cite{DudkoGrig} for a detailed account.  For $x=(x_1,x_2,\ldots)\in G^\N$, its orbit is
\[
\Gamma_Gx=\{g x\mid g\in\Gamma_G\}.
\]
The orbit can be identified with $\Gamma_G/\Stab_{\Gamma_G}(x)$.  The \emph{quasi-regular representation} $\rho_x$ on $\ell^2(\Gamma_Gx)$ is defined by
\[
(\rho_x(g)\xi)(y)=\xi(g^{-1}y).
\]
For these representations, see \cite[Sec.~2]{DudkoGrig}.
The action of $\Gamma_G$ on $G^\N$ is \emph{essentially free} when
\[
\Stab_{\Gamma_G}(x)=\{1\},\ \
\text{for}\nu\text{-almost every }x\in G^\N.
\]
For the Koopman representation $\pi$, let
\[
\mathscr H^0=L_0^2(G^\N,\nu)
=\left\{f\in L^2(G^\N,\nu)\middle|\int_{G^\N}f\,d\nu=0\right\}
=(\C\one)^\perp,\ \
\pi_0=\pi|_{\mathscr H^0},
\]
Thus $\pi_0$ is the restriction of the Koopman representation to $\mathscr H^0$, see \cite[Theorem~1(2) and Proposition~7(2)]{DudkoGrig}.

\begin{theorem}\label{t7}\label{p6}
Let $G$ be a nontrivial finite abelian group,  $\pi$ be the Koopman representation on $L^2(G^\N,\nu)$ for the action of $\Gamma_G=G\wrp\Z$ defined in Section~2, and  $\lambda$ be the left regular representation.  Then
\[
\pi_0\sim_w\pi\sim_w\lambda.
\]
\end{theorem}

\begin{proof}
The group $\Gamma_G=G\wrp\Z$ is amenable,  hence Hulanicki's theorem \cite{Hulanicki} implies $\pi\prec\lambda$, while the inclusion $\mathscr H^0\subset L^2(G^\N,\nu)$ implies that $\pi_0\prec\pi$.

The product measure $\nu$ is nonatomic.  Indeed, if $C_n(x)$ is the $n$-cylinder containing $x$, then
\[
\{x\}=\bigcap_{n\ge1}C_n(x),\ \
\nu(C_n(x))=r^{-n},
\]
so $\nu(\{x\})=0$.  By \cite[Proposition~7(2)]{DudkoGrig}, for $\nu$-almost every $x\in G^\N$, we have
\[
\rho_x\prec\pi_0.
\]
By \cite[p.~12]{KSS}, the action of $\Gamma_G$ on $G^\N$ is essentially free.
Hence there exists $x\in G^\N$ such that $\Stab_{\Gamma_G}(x)=\{1\}$ and $\rho_x\prec\pi_0$.
The orbit map $g\mapsto gx$ is a $\Gamma_G$-equivariant bijection from $\Gamma_G$ onto $\Gamma_Gx$.  Hence
\[
U_x:\ell^2(\Gamma_G)\longrightarrow\ell^2(\Gamma_Gx),\ \
U_x\delta_g=\delta_{gx},
\]
is unitary.  For all $g,h\in\Gamma_G$,
\[
U_x\lambda(h)\delta_g
=\delta_{hgx}
=\rho_x(h)U_x\delta_g,
\]
therefore $U_x\lambda(h)=\rho_x(h)U_x$ and $\rho_x\cong\lambda$.  Consequently,
\[
\lambda\prec\pi_0.
\]
Combining these containments, then
$\lambda\prec\pi_0\prec\pi\prec\lambda$ proves the asserted weak equivalences.
\end{proof}

By \cite[Sec.~2.3]{DudkoGrig}, if $\pi_1\prec\pi_2$, then $p(D_{\pi_1,G})\subseteq p(D_{\pi_2,G})$. From Theorem~\ref{t7}, $\pi\sim_w\lambda$, we have
\[
p(D_{\pi,G})=p(D_{\lambda,G}).
\]

We next determine this set explicitly.
For $u,v\in\C$, let
\[
E(u,v)=\{-r(uw+v\overline w)\mid |w|\le1\}.
\]
Since $E(cu,cv)=cE(u,v)$ for $c\in\C^\times$, the condition
$\theta\in E(u,v)$ is independent of the choice of homogeneous representative.
Assume first that $uv\ne0$ and let $\varrho_1,\varrho_2$ be the roots of
\eqref{e19}.  Let
\[
p_\theta(\xi)=ru\xi^2+\theta\xi+rv.
\]
Its zeros are $\xi_j=-\varrho_j/(ru)$.  For $|\xi|=1$,
$p_\theta(\xi)=0$ if and only if
$\theta=-r(u\xi+v\overline\xi)$.  Hence a zero of $p_\theta$ lies on the unit
circle exactly when $\theta\in\partial E(u,v)$.

If $|u|\ne|v|$, $\partial E(u,v)$ is a nondegenerate ellipse.  By the
argument principle \cite[Chap.~V, Theorem~3.4]{Conway}, the number $N(\theta)$
of zeros of $p_\theta$ in the unit disk, for
$\theta\notin\partial E(u,v)$, is
\[
N(\theta)=1+\operatorname{wind}\bigl(\theta+r(u\xi+v\xi^{-1}),0\bigr),
\ \ |\xi|=1.
\]
Thus
\[
N(\theta)=
\begin{cases}
2,& \theta\in\operatorname{int}E(u,v),\quad |u|>|v|,\\
0,& \theta\in\operatorname{int}E(u,v),\quad |u|<|v|,\\
1,& \theta\notin E(u,v).
\end{cases}
\]
If $|u|=|v|$, $E(u,v)$ is a line segment. Outside this segment there are no
zeros on the unit circle and $|\xi_1\xi_2|=|v/u|=1$, so exactly one zero lies
in the unit disk.  After ordering the roots by
$|\varrho_1|\le|\varrho_2|$, we obtain
\begin{equation}\label{e57}
\theta\notin E(u,v)
\iff
|\varrho_1|<r|u|<|\varrho_2|.
\end{equation}

\begin{theorem}\label{t8}
For every nontrivial finite abelian group $G$, one has
\[
p(D_{\pi,G})=p(D_{\lambda,G})
=\bigl\{[z]\in\Pj^3\mid\theta\in E(u,v)\bigr\}
\cup \mathcal J(\widehat F),
\]
where $\theta=z_0+(r-1)z_1$.
\end{theorem}

\begin{proof}
By Theorem~\ref{t7}, it suffices to identify
$p(D_{\pi,G})$.

Suppose first that $uv\ne0$ and
$[z]\notin\mathcal J(\widehat F)$.  
Let $\varrho_1$, $\varrho_2$ be the roots of \eqref{e38}. Without loss of generality, assume that $|\varrho_1|<|\varrho_2|$ and that $\varrho_2$ is attracting.  Then it follows from Theorem~\ref{t6} and \eqref{e50} that
\[
[z]\notin p(D_{\pi,G})
\iff
\left|\frac{u}{\varrho_2}\right|<\frac1r,\ \
\left|\frac{v}{\varrho_2}\right|<\frac1r.
\]
Since $\varrho_1\varrho_2=r^2uv$, the two inequalities are equivalent
to
\[
|\varrho_1|<r|u|<|\varrho_2|.
\]
By \eqref{e57}, this is equivalent to
$\theta\notin E(u,v)$. Therefore, for
$uv\ne0$,
\[
[z]\notin p(D_{\pi,G})
\iff
[z]\notin\mathcal J(\widehat F)
 \ \text{and}\ 
\theta\notin E(u,v).
\]

Let $v=0$ and $[z]\notin\mathcal J(\widehat F)$.  Then $\theta\ne0$ by
Theorem~\ref{t4}.  Since $M_\eta(\delta)=\theta$, from Theorem~\ref{t6} and \eqref{e50}, we have
\[
[z]\notin p(D_{\pi,G})
\iff
\left|\frac{u}{\theta}\right|<\frac1r
\iff
|\theta|>r|u|
\iff
[z]\notin\mathcal J(\widehat F)
\ \text{and}\ 
\theta\notin E(u,0). 
\] 

The case $u=0$ is obtained in the same way.
\end{proof}

\begin{remark}
The invariant subspace $\C\mathbf 1$ gives the factor $D_0(z)$.  On the other hand, the weak equivalence $\pi_0\sim_w\pi$ implies $L_G\subset p(D_{\pi_0,G})$.

For every $n\ge0$,
\[
\mathscr H_n=\C\mathbf 1_n\oplus \mathscr H_n^0,\ \
D_n(z)=D_0(z)I_{\C\mathbf 1_n}\oplus D_n^0(z),
\]
and hence
$\det D_n(z)=D_0(z)\det D_n^0(z)$.
Similarly,
\[
L^2(G^\N,\nu)=\C\mathbf 1\oplus \mathscr H^0,\ \
D_{\pi,G}(z)
=
D_0(z)I_{\C\mathbf 1}\oplus D_{\pi_0,G}(z),
\]
and hence
$p(D_{\pi,G})=L_G\cup p(D_{\pi_0,G})$.
Since $\pi_0\sim_w\pi$ by Theorem~\ref{t7}, we have
\[
p(D_{\pi_0,G})=p(D_{\pi,G}),\ \
L_G\subset p(D_{\pi_0,G}).
\]
From this, we can also see that $L_G \subset \mathcal J(\widehat F)\cup \mathcal F_\infty$.
\end{remark}

\section{Fuglede--Kadison determinants and the Li--Thom problem}\label{s4}

\subsection{Groups with a nontrivial group-ring projection}\label{ss3}

For a discrete group $K$ admitting a nontrivial projection $\mathfrak e\in \C[K]$, set
$\Gamma_K=K\wrp\Z.$
Let $\tau_K$ denote the canonical trace on $\Nalg K$, and $t$ denote the generator of the quotient $\Gamma_K/(\bigoplus_{\mathbb Z}K)\cong\mathbb Z$.  We use the Fuglede--Kadison determinant and the trace $\Tr_d\otimes\tau_{\Gamma_K}$ from Definition~\ref{d1}.

We shall use the following three properties of the Fuglede--Kadison determinant.
\begin{enumerate}[label=\textup{(\roman*)},ref=\textup{(\roman*)}]
\item\label{k1} If $T\in\Nalg\Gamma_K$ and $U_1,U_2\in\Nalg\Gamma_K$ are unitary, then
\[
\det(U_1TU_2)=\det(T).
\]
See \cite[Lemma~1(6), pp.~528--529]{FK}.

\item\label{k2} If $T\in\Nalg\Gamma_K$ has nonzero kernel, then
\[
\det(T)=0.
\]
See \cite[p.~528 and Lemma~6]{FK}.

\item\label{k3} Let $\Lambda\le\Gamma$ be countable discrete groups and let $\mathbf{A}\in M_d(\C\Lambda)$.  If
$\lambda_\Lambda(\mathbf{A})$ is injective, then
\[
\det\bigl(\lambda_{\Gamma}(\mathbf{A})\bigr)
=
\det\bigl(\lambda_\Lambda(\mathbf{A})\bigr).
\]
See \cite[Theorem~3.14(6)]{Lueck}.
\end{enumerate}

We also recall the von Neumann dimension.
\begin{definition}\cite[Sec.~1.1]{Lueck}
Let $(\mathcal M,\tau)$ be a finite von Neumann algebra.  If a closed right $\mathcal M$-submodule
$\mathscr H_{sub} \subset L^2(\mathcal M)^d$ is the range of a projection
$P\in M_d(\mathcal M)$, its von Neumann dimension is
\[
\dim_{\mathcal M}\mathscr H_{sub}=(\Tr_d\otimes\tau)(P).
\]
In particular, if $p\in\mathcal M$ is a projection,
$\dim_{\mathcal M}\bigl(pL^2(\mathcal M)\bigr)=\tau(p)$.
\end{definition}

Let $K$ be a discrete group for which $\C[K]$ contains a nontrivial projection $\mathfrak e$, and let
\[
W=\tau_K(\mathfrak e)^{-1}>1,\ \
Y_{\mathfrak e}=t^{-1}\mathfrak e.
\]
For $j\in\Z$, let $\mathfrak e_j=t^{-j}\mathfrak e t^j$.  By Lemma~\ref{l1},
$\mathfrak e_i\mathfrak e_j=\mathfrak e_j\mathfrak e_i$, $i\ne j$.
For $m\ge1$, let
\begin{equation}\label{e58}
\mathfrak p_m=(1-\mathfrak e_0)\mathfrak e_1\cdots \mathfrak e_{m-1}(1-\mathfrak e_m),\ \
\mathcal K_m=\mathfrak p_mL^2(\Nalg\Gamma_K),\ \
\omega_m^{(W)}=\frac{(W-1)^2}{W^{m+1}}
\end{equation}
so $\mathfrak p_m$ is a projection.
Let
\[
\mathfrak q_m=(1-\mathfrak e_1)\mathfrak e_2\cdots \mathfrak e_{m-1}(1-\mathfrak e_m),\ \  m\ge2 .
\]
Then
\[
\mathfrak p_m=t\mathfrak q_{m+1}t^{-1},\ \ \forall m\ge1.
\]
Recall the Jacobi matrix $\mathcal{B}_m(a;\alpha,\beta)$ from Subsection~2.3.  For $\kappa,\beta,\gamma\in\C$, let
\[
\mathfrak f_{\kappa,\beta,\gamma}=\kappa\one+\beta Y_{\mathfrak e}+\gamma Y_{\mathfrak e}^*.
\]
Dicks and Schick give a decomposition of the self-adjoint element
$Y_{\mathfrak e}+Y_{\mathfrak e}^*$.
Since $\mathfrak p_m=t\mathfrak q_{m+1}t^{-1}$, their orthogonality and completeness relations
\cite[(3.3)--(3.5), Lemmas~3.4 and~3.6, Remark~3.7]{DicksSchick} give
\begin{equation}\label{e59}
(t^j\mathcal K_m)\perp(t^{j'}\mathcal K_{m'})
\ \ \text{if }(m,j)\ne(m',j'),
\end{equation}
\begin{equation}\label{e60}
L^2(\Nalg\Gamma_K)
=
\bigoplus_{m\ge1}\;
\bigoplus_{j=0}^{m-1}t^j\mathcal K_m.
\end{equation}
For $m\ge1$, let
\[
\mathcal L_m=\bigoplus_{j=0}^{m-1}t^j\mathcal K_m.
\]
For $h\in\mathcal K_m$ and $0\le j\le m-1$, define
\[
U_m(\varepsilon_{j+1}\otimes h)=t^jh.
\]
By \eqref{e59}, the map $U_m:\C^m\otimes\mathcal K_m\to\mathcal L_m$ is unitary.    
Moreover,
\[
\dim_{\Nalg\Gamma_K}(t^j\mathcal K_m)
=\tau_{\Gamma_K}(t^j\mathfrak p_mt^{-j})
=\tau_{\Gamma_K}(\mathfrak p_m)
=\frac{(W-1)^2}{W^{m+1}}
=\omega_m^{(W)}.
\]
Let
$
U=\bigoplus_{m\ge1}U_m$.
By \eqref{e60}, we know that $U$ is unitary.

\begin{theorem}[Dicks--Schick] \cite[Lemmas~3.4 and~3.6, Remark~3.7]{DicksSchick} \label{t9}
The unitary $U$ yields a right $\Nalg\Gamma_K$-module isomorphism
\[
L^2(\Nalg\Gamma_K)\cong\bigoplus_{m\ge1}(\C^m\otimes\mathcal K_m)
\]
under which
\[
\lambda(Y_{\mathfrak e})\cong\bigoplus_{m\ge1}(J_m\otimes I_{\mathcal K_m}),\ \
\lambda(Y_{\mathfrak e}^*)\cong\bigoplus_{m\ge1}(J_m^*\otimes I_{\mathcal K_m}),
\]
and hence
\[
\lambda(\mathfrak f_{\kappa,\beta,\gamma})\cong
\bigoplus_{m\ge1}\bigl(\mathcal{B}_m(\kappa;\beta,\gamma)\otimes I_{\mathcal K_m}\bigr),\ \
\dim_{\Nalg\Gamma_K}\mathcal K_m=\omega_m^{(W)}.
\]
\end{theorem}

By \eqref{e21},
$\det \mathcal{B}_m=\mathsf T_m(\kappa,\beta\gamma)$.  Let $P_m$ denote the orthogonal
projection of $\C^m$ onto $\ker \mathcal{B}_m$.  Then
\[
\ker\lambda(\mathfrak f_{\kappa,\beta,\gamma})
\cong
\bigoplus_{m\ge1}
\bigl(\ker \mathcal{B}_m\otimes\mathcal K_m\bigr),
\]
and
\[
\begin{aligned}
	\dim_{\Nalg\Gamma_K}\bigl(\ker \mathcal{B}_m\otimes\mathcal K_m\bigr)
	&=(\Tr_m P_m)\,\dim_{\Nalg\Gamma_K}\mathcal K_m\\
	&=\bigl(\dim_\C\ker \mathcal{B}_m\bigr)\omega_m^{(W)}.
\end{aligned}
\]
Hence
\[
\dim_{\Nalg\Gamma_K}\ker \mathfrak f_{\kappa,\beta,\gamma}
=\sum\limits_{m\ge1}\omega_m^{(W)}\dim_\C\ker \mathcal{B}_m.
\]

Since $\sum\limits_{m\ge1}m\omega_m^{(W)}
=(W-1)^2\sum\limits_{m\ge1}mW^{-m-1}=1$, we obtain
\begin{equation}\label{e61}
	\log\det\bigl(\lambda(\mathfrak f_{\kappa,\beta,\gamma})\bigr)
	=\sum\limits_{m\ge1}\omega_m^{(W)}\log|\det \mathcal{B}_m|,
\end{equation}
Equation (61) is understood in the extended real sense, with the convention that $\log 0 = -\infty$. Indeed, for $\varepsilon>0$, the orthogonal decomposition of Theorem~4.2 yields
\[
\begin{aligned}
	\tau\bigl(\log(|\lambda(\mathfrak f_{\kappa,\beta,\gamma})|+\varepsilon I)\bigr)
	&=
	\sum_{m\ge1}\omega_m^{(W)}
	\operatorname{Tr}_m
	\log(|\mathcal{B}_m|+\varepsilon I_m).
\end{aligned}
\]
Since
$\|\mathcal{B}_m\|\le |\kappa|+|\beta|+|\gamma|$,
$\sum_{m\ge1}m\omega_m^{(W)}=1$,
monotone convergence, after subtracting the logarithms from a fixed upper bound, allows us to take the limit $\varepsilon \to 0$ and obtain (61), including the case where the sum equals $-\infty$.

For $x\in\R$, define
\[
\norm{x}_{\R/\Z}:=\min_{j\in\Z}|x-j|.
\]
For $0<\alpha<1$ irrational, let
\[
\alpha=[0;a_1,a_2,\ldots],
\]
and let
\begin{equation}\label{e62}
1=q_0<q_1<q_2<\cdots
\end{equation}
be the strictly increasing denominators of its continued-fraction convergents. If the first two convergents have the same denominator $1$, we discard the first one and reindex the remaining denominators.
For the continued-fraction properties can be found in \cite{Khinchin}.
\begin{lemma}\label{t10}
Let $\beta,\gamma\in\C^\times$, choose $\upsilon\in\C^\times$ with
$\upsilon^2=\beta\gamma$, and let
$\kappa=2\upsilon\cos(\pi\alpha)$ for $\alpha\in(0,1)\setminus\Q$.  Then every
$\mathcal{B}_m$ is invertible,
\[
\ker\lambda(\mathfrak f_{\kappa,\beta,\gamma})=\{0\},\ \
\det \mathcal{B}_m
=\upsilon^m\frac{\sin((m+1)\pi\alpha)}{\sin(\pi\alpha)},
\]
and,
\begin{equation}\label{e63}
\begin{aligned}
\log\det\bigl(\lambda(\mathfrak f_{\kappa,\beta,\gamma})\bigr)
={}&\log|\upsilon|-(W-1)\log|\sin(\pi\alpha)|\\
&+(W-1)^2\sum\limits_{n\ge1}W^{-n}\log|\sin(n\pi\alpha)|.
\end{aligned}
\end{equation}
For the convergent denominators $q_k$ in \eqref{e62},
\[
\det\bigl(\lambda(\mathfrak f_{\kappa,\beta,\gamma})\bigr)>0
\iff
\sum\limits_{n\ge1}W^{-n}\log\frac1{\norm{n\alpha}_{\R/\Z}}<\infty
\iff
\sum\limits_{k\ge0}W^{-q_k}\log q_{k+1}<\infty.
\]
\end{lemma}

\begin{proof}
By \eqref{e21}, \eqref{e22}, then
$\beta\gamma=\upsilon^2$, and
$\kappa=2\upsilon\cos(\pi\alpha)$,
\[
\det \mathcal{B}_m
=\mathsf T_m(\kappa,\upsilon^2)
=\upsilon^mU_m(\cos(\pi\alpha))
=\upsilon^m\frac{\sin((m+1)\pi\alpha)}{\sin(\pi\alpha)}.
\]
Since $\alpha\notin\Q$,
\[
\sin(\pi\alpha)\ne0,\ \
\sin((m+1)\pi\alpha)\ne0
\ \ \forall m\ge1.
\]
Hence
\[
\forall m\ge1,\ \ \ker \mathcal{B}_m=\{0\},
\]
and Theorem~\ref{t9} implies
\[
\ker\lambda(\mathfrak f_{\kappa,\beta,\gamma})
\cong
\bigoplus_{m\ge1}
(\ker \mathcal{B}_m\otimes\mathcal K_m)
=\{0\}.
\]

Let $h(x)=-\log|\sin(\pi x)|$, for $x\notin\Z$.
Then
\[
\log|\det \mathcal{B}_m|
=m\log|\upsilon|-h((m+1)\alpha)+h(\alpha).
\]
Since
$\omega_m^{(W)}=(W-1)^2W^{-m-1}$,
we have
\[
\sum\limits_{m\ge1}m\omega_m^{(W)}=1,\ \
\sum\limits_{m\ge1}\omega_m^{(W)}=\frac{W-1}{W}.
\]
Substitution in \eqref{e61} yields
\[
\begin{aligned}
\log\det\bigl(\lambda(\mathfrak f_{\kappa,\beta,\gamma})\bigr)
={}&\log|\upsilon|
-(W-1)^2\sum\limits_{m\ge1}W^{-m-1}h((m+1)\alpha)
+\frac{W-1}{W}h(\alpha)\\
={}&\log|\upsilon|+(W-1)h(\alpha)
-(W-1)^2\sum\limits_{n\ge1}W^{-n}h(n\alpha),
\end{aligned}
\]
which is equivalent to \eqref{e63}.

For $x\notin\Z$,
$2\norm{x}_{\R/\Z}
\le |\sin(\pi x)|
\le \pi\norm{x}_{\R/\Z}$,
hence
\[
\left|
h(x)-\log\frac1{\norm{x}_{\R/\Z}}
\right|
\le\log\pi.
\]
Since
$\sum\limits_{n\ge1}W^{-n}<\infty$,
the two nonnegative series
\[
\sum\limits_{n\ge1}W^{-n}h(n\alpha)
\ \ \text{and}\ \
\sum\limits_{n\ge1}W^{-n}
\log\frac1{\norm{n\alpha}_{\R/\Z}}
\]
converge or diverge together.  Since all other terms in
\eqref{e63} are real numbers, we have
\begin{equation}\label{e64}
\det\bigl(\lambda(\mathfrak f_{\kappa,\beta,\gamma})\bigr)>0
\iff
\sum\limits_{n\ge1}W^{-n}
\log\frac1{\norm{n\alpha}_{\R/\Z}}<\infty
\end{equation}

By \cite[Ch.~II, Secs.~6--7]{Khinchin}, we know that
\begin{equation}\label{e65}
\frac1{q_k+q_{k+1}}
<\norm{q_k\alpha}_{\R/\Z}
<\frac1{q_{k+1}},\ \
\forall k\ge0,
\end{equation}
and whenever $q_k\le n<q_{k+1}$,
\begin{equation}\label{e66}
\norm{n\alpha}_{\R/\Z}
\ge\norm{q_k\alpha}_{\R/\Z}.
\end{equation}

Assume first that
\[
\sum\limits_{n\ge1}W^{-n}
\log\frac1{\norm{n\alpha}_{\R/\Z}}<\infty.
\]
Using the terms $n=q_k$ and the second inequality in
\eqref{e65}, then
\[
\begin{aligned}
\sum\limits_{n\ge1}W^{-n}
\log\frac1{\norm{n\alpha}_{\R/\Z}}
&\ge
\sum\limits_{k\ge0}W^{-q_k}
\log\frac1{\norm{q_k\alpha}_{\R/\Z}}\\
&\ge
\sum\limits_{k\ge0}W^{-q_k}\log q_{k+1}.
\end{aligned}
\]
Hence $\sum\limits_{k\ge0}W^{-q_k}\log q_{k+1}<\infty$.

Conversely, suppose
\[
\sum\limits_{k\ge0}W^{-q_k}\log q_{k+1}<\infty.
\]
From \eqref{e65} and \eqref{e66}, for $q_k\le n<q_{k+1}$, we obtain
\[
\log\frac1{\norm{n\alpha}_{\R/\Z}}
\le\log(q_k+q_{k+1})
\le\log2+\log q_{k+1}.
\]
Therefore
\[
\begin{aligned}
&\sum\limits_{n=q_k}^{q_{k+1}-1}W^{-n}
\log\frac1{\norm{n\alpha}_{\R/\Z}}\\
&\qquad\le
(\log2+\log q_{k+1})
\sum\limits_{n=q_k}^{\infty}W^{-n}\\
&\qquad=
\frac{W}{W-1}W^{-q_k}
(\log2+\log q_{k+1}).
\end{aligned}
\]
Since $q_{k+2}\ge q_{k+1}+q_k$, thus
\[
\sum\limits_{k\ge0}W^{-q_k}<\infty.
\]
Therefore
\[
\begin{aligned}
&\sum\limits_{n\ge1}W^{-n}
\log\frac1{\norm{n\alpha}_{\R/\Z}}\\
&\quad=
\sum\limits_{k\ge0}
\sum\limits_{n=q_k}^{q_{k+1}-1}W^{-n}
\log\frac1{\norm{n\alpha}_{\R/\Z}}\\
&\quad\le
\frac{W}{W-1}
\sum\limits_{k\ge0}W^{-q_k}
\bigl(\log2+\log q_{k+1}\bigr)<\infty.
\end{aligned}
\]
Together with \eqref{e64}, this proves the equivalences.
\end{proof}

\begin{theorem}\label{t11}
If $K$ is a discrete group for which $\C[K]$ contains a nontrivial projection, then there exists $\mathfrak f\in\C[K\wrp\Z]$ such that
\[
\ker\lambda(\mathfrak f)=\{0\},\ \
\det(\lambda(\mathfrak f))=0.
\]
\end{theorem}
\begin{proof}
Let $W=\tau_K(\mathfrak e)^{-1}>1$.  Choose $a_1\ge2$ and then select an irrational $\alpha=[0;a_1,a_2,\ldots]$ recursively so that its convergent denominators satisfy
\[
q_{k+1}=a_{k+1}q_k+q_{k-1}\ge\exp(W^{q_k})
\ \ \forall k\ge1.
\]
Then $W^{-q_k}\log q_{k+1}\ge1$, and hence $\sum\limits_{k\ge0}W^{-q_k}\log q_{k+1}=\infty$.

Let $\beta=\gamma=\upsilon=1$, $\kappa=2\cos(\pi\alpha)$, and $\mathfrak f=2\cos(\pi\alpha)\one+Y_{\mathfrak e}+Y_{\mathfrak e}^*$.  By Lemma~\ref{t10},
\[
\det \mathcal{B}_m=\frac{\sin((m+1)\pi\alpha)}{\sin(\pi\alpha)}\ne0
\]
Thus
\[
\det \mathcal{B}_m\ne0,\ \ \forall m\ge1,
\]
so $\ker\lambda(\mathfrak f)=\{0\}$.  Since $\sum\limits_{k\ge0}W^{-q_k}\log q_{k+1}=\infty$, Lemma~\ref{t10} shows that
\[
\log\det(\lambda(\mathfrak f))=-\infty,\ \
\det(\lambda(\mathfrak f))=0.
\]

If $k\in K$ has finite order $s\ge2$, then $e_k=s^{-1}\sum\limits_{j=0}^{s-1}k^j\in\C[K]$ is a nontrivial projection with $\tau_K(e_k)=1/s$.  Hence the construction applies with $\mathfrak e=e_k$ and $W=s$.
\end{proof}

\begin{example}
Let $K=\Z_2=\langle a\mid a^2=1\rangle$ and $e=(\one+a)/2$, so $W=2$. Let
\[
N_1=1,\ \ N_2=5,\ \
N_{j+1}=N_j+\left\lceil 2^{10^{N_j}}\right\rceil, \ \ j\ge2,
\]
and
\[
\alpha_*=\sum_{j\ge1}10^{-N_j}
=0.10001\underbrace{00\cdots0}_{N_3-N_2-1}1\cdots .
\]
Since $N_{j+1}-N_j\to\infty$, $\alpha_*\notin\Q$.  For $j\ge2$,
\[
\norm{10^{N_j}\alpha_*}_{\R/\Z}
\le 2\times10^{N_j-N_{j+1}},
\]
and therefore
\[
2^{-10^{N_j}}
\log\frac1{\norm{10^{N_j}\alpha_*}_{\R/\Z}}
\ge
\log 10-2^{-10^{N_j}}\log2.
\]
Hence the series in \eqref{e64} diverges.  Therefore, for
\[
\mathfrak f=2\cos(\pi\alpha_*)\one
+\frac12t^{-1}(\one+a)+\frac12(\one+a)t
\in\C[\Z_2\wrp\Z],
\]
Lemma~\ref{t10} implies that
\[
\ker\lambda(\mathfrak f)=\{0\},\ \
\det(\lambda(\mathfrak f))=0.
\]
\end{example}

\begin{remark}
Theorem~\ref{t11} implies scalar elements satisfying
\[
\ker\lambda(\mathfrak f)=\{0\},\ \
\det(\lambda(\mathfrak f))=0.
\]
If $K$ is countable amenable, then $K\wrp\Z$ is amenable, and Theorem~\ref{t11} gives a negative answer to Li--Thom Question~3.13 for these wreath products.  In particular, every nontrivial finite $K$ is covered.
\end{remark}

\subsection{Free abelian base groups: injectivity implies positive determinant}\label{ss4}

Fix $m\ge1$ and let
\[
\Gamma_m=\Z^m\wrp\Z=\mathcal{G}_m\rtimes\langle t\rangle,\ \
\mathcal{G}_m=\bigoplus_{j\in\Z}(\Z^m)_j,
\]
where $\Z^m$ denotes the free abelian group of rank $m$.
Every element of $\Gamma_m$ has a unique form $bt^k$ with $b\in \mathcal{G}_m$ and $k\in\Z$, see \cite[Chap.~7]{RotmanGroups}.  Hence every nonzero $\mathfrak a\in\C\Gamma_m$ has a unique finite expansion
\[
\mathfrak a=\sum_{k=j}^{n}c_kt^k,\ \
c_k\in\C \mathcal{G}_m,\ \
c_j\ne0.
\]
Let $\widehat{\mathcal{G}_m}$ be the Pontryagin dual of $\mathcal{G}_m$ and let $\mu_m$ be normalized Haar measure on $\widehat{\mathcal{G}_m}$.  The Fourier transform $\mathcal F_{\mathcal{G}_m}$ shows that
\begin{equation}\label{e67}
\mathcal F_{\mathcal{G}_m}\,\Nalg \mathcal{G}_m\,\mathcal F_{\mathcal{G}_m}^{-1}
=L^\infty(\widehat{\mathcal{G}_m},\mu_m).
\end{equation}
See \cite[Chap.~1]{RudinFA}.  
By \eqref{e67} and \cite[Chap.~X]{TakesakiII}, we have
\begin{equation}\label{e68}
	\Nalg\Gamma_m\cong L^\infty(\widehat{\mathcal{G}_m},\mu_m)\rtimes\Z.
\end{equation}


Let $c\in\C \mathcal{G}_m$ be nonzero.  Then there exist an integer $M\ge1$, pairwise distinct elements $b_1,\ldots,b_M\in \mathcal{G}_m$, and coefficients $c_1,\ldots,c_M\in\C\setminus\{0\}$ such that
\[
c=\sum_{j=1}^{M}c_jb_j.
\]
Let
$\Lambda=\langle b_1,\ldots,b_M\rangle$.
Then $\Lambda\cong\Z^N$ for some $N\ge0$.  Choose a basis $e_1,\ldots,e_N$ of $\Lambda$.  For $1\le j\le M$ there are unique integers $m_{1j},\ldots,m_{Nj}$ such that
\[
b_j=e_1^{m_{1j}}\cdots e_N^{m_{Nj}}.
\]
Let
\[
\mathbb T=\{z\in\C\mid |z|=1\}.
\]
For $\chi\in\widehat\Lambda$, let $z_i=\chi(e_i)$.  Then
\[
\widehat\Lambda\cong\mathbb T^N,\ \
\chi\longmapsto (z_1,\ldots,z_N).
\]
Under the Fourier transform, $\lambda_\Lambda(c)$ is multiplication by the Laurent polynomial
\[
P_c(z_1,\ldots,z_N)
=\sum_{j=1}^{M}c_jz_1^{m_{1j}}\cdots z_N^{m_{Nj}},\ \
(z_1,\ldots,z_N)\in\mathbb T^N.
\]
Let $\mu_N$ be normalized Haar measure on $\mathbb T^N$.  The Mahler-measure formula implies that
\[
\log|P_c|\in L^1(\mathbb T^N,\mu_N),\ \
\int_{\mathbb T^N}\log|P_c(z)|\,d\mu_N(z)>-\infty.
\]
See \cite[pp.~45--46]{DeningerMahler}.  Therefore
\[
P_c(z)\ne0
\ \ \text{for $\mu_N$-almost every }z\in\mathbb T^N,
\]
and multiplication by $P_c$ is injective.  Hence
\[
\ker\lambda_\Lambda(c)=\{0\},\ \
\log\det\bigl(\lambda_\Lambda(c)\bigr)
=
\int_{\mathbb T^N}\log|P_c(z)|\,d\mu_N(z)>-\infty.
\]
By property~\ref{k3} in Subsection~\ref{ss3},
\[
\log\det\bigl(\lambda_{\mathcal{G}_m}(c)\bigr)
=
\log\det\bigl(\lambda_\Lambda(c)\bigr)>-\infty.
\]
From property~\ref{k2} in Subsection~\ref{ss3}, we have
\begin{equation}\label{e69}
\ker\lambda_{\mathcal{G}_m}(c)=\{0\},\ \
\det\bigl(\lambda_{\mathcal{G}_m}(c)\bigr)>0.
\end{equation}

\begin{theorem}\label{t12}
Let
\[
 \mathfrak a=\sum_{k=j}^{n}c_kt^k\in\C\Gamma_m
\]
be nonzero, with $c_j\ne0$.  Then
\[
\ker\lambda(\mathfrak a)=\{0\},\ \
\log\det\bigl(\lambda_{\Gamma_m}(\mathfrak a)\bigr)
\ge
\log\det\bigl(\lambda_{\mathcal{G}_m}(c_j)\bigr)>-\infty.
\]
In particular,
\[
\det\bigl(\lambda_{\Gamma_m}(\mathfrak a)\bigr)>0.
\]
\end{theorem}

\begin{proof}
Applying \eqref{e69} to $c_j$ yields
$\ker\lambda_{\mathcal{G}_m}(c_j)=\{0\}$, $\det\bigl(\lambda_{\mathcal{G}_m}(c_j)\bigr)>0$.
Let
\[
\mathfrak b=\mathfrak at^{-j}
=c_j+c_{j+1}t+\cdots+c_nt^{n-j}.
\]
Since $\lambda(t)$ is unitary, property~\ref{k1} in Subsection~\ref{ss3} implies
\[
\det\bigl(\lambda_{\Gamma_m}(\mathfrak b)\bigr)
=
\det\bigl(\lambda_{\Gamma_m}(\mathfrak a)\bigr).
\]

By \eqref{e68}, let
\[
\mathscr{A}_+
=
\overline{
\left\{
\sum_{k=0}^{N}a_kt^k
\ \middle|\ 
N\ge0,\ a_k\in L^\infty(\widehat{\mathcal{G}_m},\mu_m)
\right\}
}^{\,w^*}
\subset\Nalg\Gamma_m,
\]
and let
\[
\mathbb E:\Nalg\Gamma_m\longrightarrow\Nalg \mathcal{G}_m,\ \
\mathbb E\left(\sum_{k\in\Z}a_kt^k\right)=a_0
\]
be the trace-preserving conditional expectation.  With respect to $\mathbb E$, the algebra $\mathscr{A}_+$ is a finite subdiagonal algebra of the crossed product, and it is maximal by Exel's theorem \cite[Theorem~7]{Exel}. By Jensen's inequality
for finite maximal subdiagonal algebras  \cite[p.~3644 and Theorem~3]{Labuschagne}, we have
\[
\det\bigl(\lambda_{\Gamma_m}(\mathfrak b)\bigr)
\ge
\det\bigl(\mathbb E\bigl(\lambda_{\Gamma_m}(\mathfrak b)\bigr)\bigr)
=
\det\bigl(\lambda_{\mathcal{G}_m}(c_j)\bigr)>0.
\]
Then, 
$\det\bigl(\lambda_{\Gamma_m}(\mathfrak a)\bigr)>0$.
By property~\ref{k2} in Subsection~\ref{ss3},  therefore
$\ker\lambda(\mathfrak a)=\{0\}$.
\end{proof}

We next consider matrices.  Let
$S=\C\Gamma_m\setminus\{0\}$.
Since $\Gamma_m$ is torsionfree elementary amenable, $\C\Gamma_m$ has no nontrivial zero divisors, and its skew field of fractions is $S^{-1}\C\Gamma_m$, obtained via Ore localization with respect to $S$ \cite[Remark~5.1]{LueckTwist}.
For $d\ge1$ and $\mathbf{A}\in M_d(\C\Gamma_m)$, let
\begin{equation}\label{e70}
r_\mathbf{A}:\ell^2(\Gamma_m)^d\longrightarrow\ell^2(\Gamma_m)^d
\end{equation}
denote right multiplication by $\mathbf{A}$ \cite[Sec.~7]{LueckTwist}.
For $x\in\C\Gamma_m$ with $x\ne0$, the operator $r_x$ is said to be of determinant class if
\begin{equation}\label{e71}
\int_{(0,1)} |\log s|\,d\mu_{|r_x|}(s)<\infty.
\end{equation}
Equivalently, $\det r_x>0$, which can be found in \cite[Sec.~2]{LueckTwist}.

\begin{theorem}\label{t13}
If $\mathbf{A}\in M_d(\C\Gamma_m)$ and $\ker\lambda(\mathbf{A})=\{0\}$, then
\[
\det\lambda(\mathbf{A})>0.
\]
\end{theorem}

\begin{proof}
Let
\[
(\mathscr J\xi)(g)=\overline{\xi(g^{-1})},\ \
\mathscr J_d=\mathscr J^{\oplus d}.
\]
Then $\mathscr J$ and $\mathscr J_d$ are antiunitary and, with the notation in \eqref{e70},
\[
\mathscr J\lambda(x)\mathscr J^{-1}=r_{x^*},\ \
\mathscr J_d\lambda(A)\mathscr J_d^{-1}=r_{A^*},
\]
where $\mathbf{A}^*=(a_{ji}^*)$.
Hence
\begin{equation}\label{e72}
\ker r_{\mathbf{A}^*}=\{0\}.
\end{equation}

Let $x\in\C\Gamma_m$ with $x\ne0$.  Applying Theorem~\ref{t12} to $x^*$ shows that 
\[
\ker r_x=\{0\},\ \
\det r_x=\det\lambda(x^*)>0.
\]
Thus Assumption~7.13 of \cite{LueckTwist} follows from \eqref{e71}.

We apply \cite[Lemma~7.23]{LueckTwist} with field $\C$ and the one-dimensional trivial representation.  By \eqref{e72} and Lemma~7.23(1), then $\mathbf{A}^*$ is invertible over $S^{-1}\C\Gamma_m$. From Lemma~7.23(2), there exist $x,y\in\C\Gamma_m$, $x\ne0$, $y\ne0$, such that
\[
\det r_{\mathbf{A}^*}
=
\frac{\det r_x}{\det r_y}>0.
\]
Since antiunitary conjugation preserves singular values, then
$\det\lambda(\mathbf{A})=\det r_{\mathbf{A}^*}>0$.
\end{proof}

More generally, Theorems~\ref{t12} and \ref{t13} carry over to
$\Gamma=G_{tf}\rtimes_\sigma\mathbb Z$,
where $G_{tf}$ is any countable torsion-free abelian group and
$\sigma\in\operatorname{Aut}(G_{tf})$. Indeed, since every finitely generated subgroup of $G_{\mathrm{tf}}$ is free abelian, the Mahler-measure and Jensen-inequality arguments remain valid. Furthermore, $\Gamma$ is torsion-free elementary amenable, and the matrix argument based on \cite[Lemma~7.23]{LueckTwist} still applies. Consequently, we obtain the following corollary.

\begin{corollary}\label{c14}
Let $\Gamma=G_{tf}\rtimes_\sigma\mathbb Z$, where $G_{tf}$ is any countable torsion-free
abelian group and $\sigma\in\operatorname{Aut}(G_{tf})$. For every $d\ge1$ and every $\mathbf{A}\in M_d(\C\Gamma)$, injectivity of $\lambda(\mathbf{A})$ implies $\det\lambda(\mathbf{A})>0$.
\end{corollary}

\begin{remark}
Corollary~\ref{c14} gives an affirmative answer to Li and Thom Question~3.13 for $\Gamma=G_{tf}\rtimes_\sigma\mathbb Z$ in every matrix size.  In contrast, Theorem~\ref{t11} gives scalar counterexamples for wreath products $K\wrp\Z$ when $\C[K]$ contains a nontrivial projection.  The two theorems concern different classes of semidirect products and do not give a general dichotomy between torsion and torsion-free groups.
\end{remark}

\section{Concluding Remark}

We conclude this paper with two directions for future research. The first concerns the generality of the inclusion of the Julia set in the projective spectrum. In the present setting, this inclusion follows from the
density of the singular hypersurfaces $\mathcal{H}_n$ of the finite pencils
$D_n(z)$ in the Julia set. That mechanism comes from the Schur recursion and
the self-similar structure of the action, not from special features of
$G \wr \mathbb{Z}$. This leads us to the following conjecture.
 
\begin{conjecture}
Let $\pi$ be a self-similar representation of a wreath product group, and let
$F$ be the rational map obtained from the pencil $D_\pi$ by Schur reduction.
Then
\[
  \mathcal{J}(F) \subseteq p(D_\pi).
\]
\end{conjecture}
 
Natural test cases include wreath products $G \wr H$ with $G$ a non-abelian
finite group or with $H \neq \mathbb{Z}$, as well as other self-similar
realizations of these groups. Our results also suggest that the reverse
inclusion should fail in general for non-normal pencils. A finer problem is
then to determine which Fatou components meet the spectrum, and whether they
can be read off from the multipliers and normalized pencil coefficients
at the attracting fixed points, as happens
for the set $\mathcal{F}_\infty$ here.

The second direction concerns the Li--Thom question. By Theorem~\ref{t11}, the answer
is negative for $K\wr\mathbb Z$ whenever $\mathbb C[K]$ contains a nontrivial projection, whereas Corollary~\ref{c14} gives a positive answer for $G_{tf}\rtimes_\sigma\mathbb Z$ with $G_{tf}$ countable torsion-free abelian. This raises the following question.

\begin{question}
For which groups $G$ and automorphisms $\sigma \in \operatorname{Aut}(G)$ does every
$A \in M_d\bigl(\mathbb{C}[G \rtimes_\sigma \mathbb{Z}]\bigr)$ with
$\ker \lambda(A) = \{0\}$ satisfy $\det \lambda(A) > 0$?
\end{question}
 
A positive answer for $G \rtimes_\sigma \mathbb{Z}$ forces one for $G$ itself,
since the determinant does not change under passage to the subgroup $G$. The question is which cyclic extensions preserve this property and how the answer depends on $\sigma$. In particular, one may ask whether, for amenable $G$, the absence of nontrivial
projections in $\mathbb{C}[G]$ is sufficient. Torsion-free elementary amenable groups, for which Ore localization is available, seem a natural first class to test.

\section*{Acknowledgements}
The first-named author is supported by the National Natural Science Foundation of China (NOs.12271298 and 12671156) and the Taishan Scholars Project (TSTP 20250702). The third-named author is supported in part by the Simons Foundation Grant No. MP-TSM-00002315, the Institute for Advanced Study Summer Collaboration Program (2025), and the American Institute of Mathematics SQuaRE Program (2025-2027).

\vspace{2em}
\begin{flushleft}
	\textbf{Penghui Wang}\\
	\textit{School of Mathematics, Shandong University, Jinan 250100, China}\\
	Email address: \href{mailto:phwang@sdu.edu.cn}{phwang@sdu.edu.cn}
	
	\vspace{1em}
	
	\textbf{Xumin Wang}\\
	\textit{School of Mathematics, Shandong University, Jinan 250100, China}\\
	Email address: \href{mailto:xm.wang@mail.sdu.edu.cn}{xm.wang@mail.sdu.edu.cn}
	
	\vspace{1em}
	
	\textbf{Rongwei Yang}\\
	\textit{Department of Mathematics and Statistics, SUNY at Albany, Albany, NY 12047, USA}\\
	Email address: \href{mailto:ryang@albany.edu}{ryang@albany.edu}
\end{flushleft}

\end{document}